\documentstyle{amsppt}

\catcode`\@=11
\def\proclaimheadfont@{\smc}
\def\demoheadfont@{\smc}
\catcode`\@=\active
\expandafter\redefine\csname logo\string @\endcsname{}

\magnification=\magstep1 \pagewidth{16 true cm} \pageheight{23,5 true cm}



\topmatter

\def \Om{\Omega}
\def \ti{\times}
\def \sb{\subset}
\def \ov{\overline}
\define \>{\geqslant}
\define \<{\leqslant}

\def \al{\alpha}
\def \lm{\lambda}
\def \pri{\prime}
\def \ve{\varepsilon}
\def \pr{\partial}
\def \tl{\tilde}
\def \g{\gamma}
\mathchardef\emptyset="001F

\title
A MONOTONICITY APPROACH TO NONLINEAR DIRICHLET
PROBLEMS IN PERFORATED DOMAINS
\endtitle

\leftheadtext{Gianni Dal Maso, Igor V.~Skrypnik}
\rightheadtext{Nonlinear Dirichlet Problems in perforated domains}

\author
Gianni Dal Maso ($^{*}$), Igor V.~Skrypnik ($^{**}$)
\endauthor

\abstract
We study the asymptotic behaviour of solutions to Dirichlet problems 
in perforated domains for nonlinear elliptic equations associated with 
monotone operators. The main difference with respect to the previous 
papers on this subject is that no uniformity is assumed in the 
monotonicity condition. Under a very general hypothesis on the holes
of the domains, we construct a limit equation, which is satisfied 
by the weak limits of the solutions. The additional term 
in the limit problem depends only on the local behaviour of the holes,
which can be expressed in terms of suitable nonlinear capacities associated 
with the monotone operator. 
\endabstract

\affil
($^{*}$) SISSA, via Beirut 4, 34014 Trieste, Italy\\
e-mail: {\tt dalmaso\@sissa.it}\\
($^{**}$) Institute of Applied Mathematics and Mechanics,\\
Academy of Sciences of Ukraine, R. Luxemburg St. 74,\\
340114 Donetsk, Ukraine\\
e-mail: {\tt skrypnik\@iamm.ac.donetsk.ua}
\endaffil

\endtopmatter

\document

\baselineskip=16pt
\TagsOnLeft

\head {Introduction}
\endhead

This paper continues previous investigations of the authors on nonlinear
Dirichlet problems in perforated domains of general structure.

Let $\Om$ be any bounded domain in the  $n$-dimensional Euclidean
space $R^n$ and let $\Om_s\sb\Om,\,\,\,s=1,2,...$ be a sequence of
subdomains. In $\Om_s$ we consider a nonlinear elliptic boundary
value problems for $\,\,s=1,2,...$
$$
\underset{j=1}\to{\overset{n}\to\sum} \,\,\frac d{d x_j} \,a_j\,
\big (x, \frac{\pr u}{\pr x}\big )\, = \,
\underset{j=1}\to{\overset{n}\to\sum} \,\,\frac {\pr}{\pr x_j}\,\,
f_j (x), \,\,\,x\in\Om_s,
\tag{0.1}
$$
\vskip-10pt
$$
u(x) = f(x),\,\,\,x\in\pr\Om_s.
\tag{0.2}
$$

Our conditions on the data of problems (0.1), (0.2) provide the
existence of a solution $u_s(x) \in W_m^1 (\Om_s)$ for every $\,\,s\,\,$ 
and also the boundedness of the
sequence $u_s(x)\,\,$ in $\,\,W_m^1 (\Om)$.

In the previous works on this subject [1-4, 6-10] 
(see also the References in [6]) the homogenization problems
for nonlinear elliptic second order equations were studied under strong
monotonicity assumption for the equations. The following inequality
$$
\underset{j=1}\to{\overset{n}\to\sum} \,\,[a_j (x,p) - a_j (x,q)]
\,\,(p_j - q_j) \> \nu \,|p - q|^m
\tag{0.3}
$$
was assumed for arbitrary $p, q \in R^n,\,\,x \in \Om$, with a positive
constant $\,\nu$.

In particular, inequality (0.3) guarantees in [6-10] the strong
convergence to zero in $W_m^1 (\Om)$ of the remainder term of the
asymptotic expansion and the strong convergence of the gradients of
solutions of the problems (0.1), (0.2) in $W_p^1 (\Om)$ for
$p < m$.

In this paper we assume only the following weak monotonicity condition:
for arbitrary points $x \in \Om,\,\, p, q \in R^n$ the inequality
$$
\underset{j=1}\to{\overset{n}\to\sum} \,\,[a_j (x,p) - a_j (x,q)]\,\,
(p_j - q_j) \> 0
\tag{0.4}
$$
is satisfied. This weak condition does not allow us to apply the
methods from [6-10], which are based on the study of the behaviour of 
the asymptotic expansion of the solutions. We develop a new approach by 
using monotonicity arguments. This allows us to construct a boundary 
value problem (in a fixed domain), which is satisfied
by the weak limits of subsequences of $u_s (x)$.

This approach is based on the construction of special test functions
and on the analysis of their behaviour. For this
analysis we use precise pointwise and integral estimates for the
potential functions which are solutions of some auxiliary boundary 
value problems in domains with holes of small diameter.

We call the attention of the reader to the main result of this 
paper, the Convergence Theorem (Theorem 1.1), that is proved by
using a new pointwise estimate (Lemma 2.3) of the potential functions.
This theorem allows us to make the main modification 
in the construction of the corrector, if assumption (0.3) is satisfied. 
Note that in the previous papers [6-10] the
definition of the subdivision of the domain, and consequently the
construction of the corrector, depended on sequence $u_s (x)$.
The subdivision and the corrector we construct in the present paper, by 
using the Convergence Theorem, are independent of $u_s (x)$.

Our assumption on the perforated domains (see Condition B in
Section 1) coincides with the corresponding condition in [10]. We 
suppose that the $C_m$-capacity of the portion of the holes in any small cube 
is estimated from above by the Lebesgue measure of the cube.

We construct the limit boundary value problem, 
and we describe the additional term which appears 
in it by means of some
quantitative capacitary properties of the holes.


\head
{1. Statement of the results}
\endhead

We assume that the functions $\,\,a_j (x,p),\,\, j = 1,...,n\,$,
are defined for $x \in R^n$, $p \in R^n$, and satisfy the
following conditions:
\medskip
{\bf Condition A.1.} {\it The functions $\,\,a_j (x,p)\,\,$
are continuous in $\,\,p\,\,$ for all $\,\,x \in R^n\,\,$ and
measurable in $\,\,\,x\,\,\,$ for all $\,\,p \in R^n$.}
\medskip
{\bf Condition A.2.} {\it There exist two positive constants
$\,\,\,\nu_1, \nu_2\,$, and a constant $m$, with $\,\,2\< m < n$,
such that
$$
\underset{j=1}\to{\overset{n}\to\sum} \,a_j\,(x, p)\,p_j \>
\nu_1 \,(1 + |p|)^{m-2} \cdot |p|^2,
\tag {1.1}
$$
\vskip-20pt
$$
\underset{j=1}\to{\overset{n}\to\sum} \, \big[ a_j\,(x,p) -
a_j\,(x,q)\big]\,\,(p_j -q_j) \> 0,
\tag{1.2}
$$
\vskip-20pt
$$
\underset{j=1}\to{\overset{n}\to\sum} \,\big|a_j\,(x,p) -
a_j\,(x,q)\big| \< \,\nu_2 \,(1 + |p| + |q|)^{m-2} \,\cdot\,|p - q|,
\tag{1.3}
$$
for every $x \in R^n,\,\,\,p, q \in R^n\,\,$.}
\medskip

Note that from (1.1) it follows that $a_j\,(x,0)=0$ for 
every $x\in R^n\,\,$. Therefore (1.3) implies that
$$
|a_j\,(x,p)| \< \,\nu_2 \,(1 + |p|)^{m-2} \,|p|
\tag{1.4}
$$
for every $\,\,x \in R^n, \,\,p \in R^n, \,\,j = 1,...,n$.

We assume that functions $\,\,f_j (x),\,\,j=1,...,n\,$, and $\,\,f(x)\,\,$
in (0.1), (0.2) are defined in $\,\,R^n\,\,$ and satisfy the conditions:
$$
f_j (x) \in L_{m^\pri} (R^n),\,\,\quad f (x) \in W_m^1 (R^n)
\tag{1.5}
$$
for $\,\,j = 1,...,n\,\,$ and $\,\,m^\pri = \frac m{m-1}$.

A solution of the boundary value problem (0.1), (0.2) is a
function $u (x) \in W_m^1 (\Om _s)$, satisfying
 $\,\,u (x) - f (x) \in \overset\circ\to W_m^1 (\Om _s)\,$,
such that the integral identity
$$
\underset{j=1}\to{\overset{n}\to\sum} \,\,
\underset{\Om_s}\to\int \,\,\bigg[a_j\,\big(x , \frac {\pr u}{\pr x}\big)
- f_j (x)\bigg]\,\,\frac {\pr\varphi (x)}{\pr x_j}\,\,dx = 0
\tag{1.6}
$$
holds for an arbitrary function
$\varphi (x) \in \overset\circ\to W_m^1 (\Om _s)$.

Using methods of the theory of monotone operators it is easy to prove
the existence of a solution of problem (0.1), (0.2). For every
$\,\,s\,\,$ we denote by $\,u_s (x)\,$ one of the possible solution of
the problem (0.1), (0.2) and extend $\,\,u_s (x)\,\,$ on $\,\,R^n\,\,$
by setting $\,\,u_s (x) = f(x)\,\,$ for $\,\,x \in R^n \setminus \Om_s$.
By condition A.2 and (1.5) the estimate
$$
\underset {R^n}\to\int\,\,\bigg\{\big|\frac{\pr u_s (x)}{\pr x}\big|^m
\, + \, |u_s (x)|^m\,\bigg\}\,\,dx \< R
\tag{1.7}
$$
holds with a constant $\,\,R\,\,$ independent of $\,\,s\,$.

By (1.7) the sequence $\,u_s (x)\,$ contains a weakly convergent
subsequence, therefore we may assume that $\,u_s (x)\,$ converges weakly in
$\,W_m^1 (R^n)\,$ to some function $u_0 (x)$.

We formulate now our assumptions on the sequence $\Om_s$ in terms of the
$m$-capacity $C_m (F)$. For every compact set $\,\,F\,$, its
$m$-capacity $C_m (F)$ is defined by
$$
C_m (F) = \inf\,\underset {R^n}\to\int
\big|\frac{\pr \varphi (x)}{\pr x}\big|^m\,\,dx\,,
\tag{1.8}
$$
where the infimum is taken over all function
$\varphi (x) \in C_0^\infty (R^n)$ which satisfy the condition
$\varphi (x) = 1\,\,$ for $\,\,x \in F$.

For every $x_0 = (x_1^{(0)},...,x_n^{(0)}) \in R^n$, $r > 0$, we set
$$
K (x_0, r) = \{x \in R^n : |x_j - x_j^{(0)}| \< r, \,\,j = 1,...,n\}.
\tag{1.9}
$$
Let us assume that the following condition is satisfied.
\medskip
{\bf Condition B.} {\it There exist a positive number $A$ and a sequence
$r_s>0$, tending to zero as $s\to\infty$, such that the inequality
$$
C_m (K (x,r)\setminus \Om_s) \< A\,r^n
\tag{1.10}
$$
holds for every $x \in \Om$ and for every $r \> r_s$ with
$K (x, r+r_s) \sb\Om$.}
\medskip

Let us fix a bounded open set $\Om_0 \sb R^n$ such that
$\rho (\pr \Om_0, \Om) \> 1$, where $\rho (\pr \Om_0, \Om)$ is the distance
from $\pr \Om_0\,\,$ to $\,\Om\,$, and let $\psi (x)$ be a function
of class $C_0^\infty (\Om_0)$ equal to 1 on $\ov\Om$. For every
compact set $F$ contained in $\Om$ and for every real
number $\,\,q\,\,$ we define the auxiliary function
$v (x,F,q)$ as a solution of the boundary value problem
$$
\underset{j=1}\to{\overset{n}\to\sum}\,\frac\pr{\pr x_j}\,\,
a_j\big(x, \frac {\pr v}{\pr x}\big) = 0, \quad x \in \Om_0\setminus 
F\,,
\tag{1.11}
$$
\vskip-20pt
$$
v (x) = q\,\psi (x), \quad x \in \pr (\Om_0\setminus F)\,.
\tag{1.12}
$$
The solvability of problem (1.11), (1.12) follows easily from the
theory of monotone operators. In [5] it is proved that this problem 
admits a maximal solution, i.e., there exists a
solution $\ov v (x)$ of problem (1.11), (1.12)
such that $v (x) \< \ov v (x)$ for any solution $v (x)$ of the same
problem. We denote this maximal solution by
$v (x, F, q)$, and extend it to $\,\,R^n$ by setting
$v (x, F, q) = q$ in $F$ and $v (x, F, q) = 0$ outside $\Om_0$.

In Section 3 we shall introduce a special decomposition of the domain
$\Om$ of the form
$$
\Om = \big\{ \underset {\al\in I_s}\to\cup\,\,K (x_\al^{(s)},
\lm_s \rho_s)\big\} \cup U_s
\tag{1.13}
$$
where $\lm_s\,\,$ and $\,\,\rho_s\,$ are sequences of positive real numbers
such that $\lm_s \to \infty,\,\,\rho_s\to 0\,\,$ and
$\,\lm_s\rho_s\to 0\,\,$ as $\,\,s\to\infty,\,\,
x_\al^{(s)} = 2\lm_s \rho_s \al,\,\, \al= (\al_1,...,\al_n)$ is a
multi-index with integer coordinates, $I_s$ is the set of all
multi-indices $\,\,\al\,\,$ such that $K (x_\al^{(s)}, 2\lm_s \rho_s)\sb\Om$,
and $\,\,U_s\,\,$ is the complement of
$ \underset {\al\in I_s}\to\cup\,\,K (x_\al^{(s)}, 2\lm_s \rho_s)$
with respect to $\Om$.

We define $v_\al^{(s)} (x,q) = v (x,F,q)\,\,$ for
$\,\,F = K (x_\al^{(s)}, (\lm_s - 2)\rho_s) \setminus \Om_s$. Let
$q_s (x)$ be an arbitrary sequence that converges strongly in
$L_m (\Om)$ and let $q_\al^{(s)}$ be the mean value of the function
$q_s (x)$ in the cube $K(x_\al^{(s)}, \lm_s \rho_s)$.

In Section 3 we shall construct the following sequence, which is fundamental
in our analysis:
$$
r_s (x) = \underset {\al\in I_s}\to\sum\,\,
v_\al^{(s)} (x, q_\al^{(s)})\,\varphi_\al^{(s)} (x)\,,
\tag{1.14}
$$
where $\varphi_\al^{(s)} (x)$ is a special cut-off function,
constructed by using $v_\al^{(s)} (x, q_\al^{(s)})$ (see (3.6)),
which is equal to 1 for
$x \in K (x_\al^{(s)}, (\lm_s - 2)\rho_s) \setminus \Om_s$ and equal to
$\,\,0\,\,$ outsite  $K(x_\al^{(s)}, \lm_s \rho_s)$. Remark that
$r_s (x)$ is analogous with the corrector which was constructed in [6,9].
In Section 3 we shall prove the following result.
\medskip
{\bf Theorem 1.1 (Convergence Theorem).} {\it Assume that conditions
A.1, A.2, and B are satisfied and let $q_s (x)$ be some sequence converging
strongly in $L_m (\Om)$. Let $z_s (x)$ be an arbitrary sequence of
functions such that $z_s (x) \in \overset\circ\to W_m^1 (\Om_s)$ and
$z_s (x)$ converges weakly to zero in $W_m^1 (\Om)$. Then
$$
\underset {s\to\infty}\to\lim\,\,\underset{j=1}\to{\overset{n}\to\sum}
\,\,\underset\Om\to\int\,\,a_j\big(
x, \frac {\pr r_s (x)}{\pr x}\big) \,\,\frac{\pr z_s (x)}{\pr x_j}\,\,dx = 0.
\tag{1.15}
$$
}
\medskip

In order to formulate a result about the boundary value problem for the
function $u_0 (x)$ we introduce a capacity connected with the
differential equation (0.1), defined for every compact set $F\sb\Om$
and for every real number $q \ne 0$ by the equality
$$
C_A (F,q) = \underset{j=1}\to{\overset{n}\to\sum}\,\,\frac 1q\,\,
\underset\Om\to\int\,\,a_j\big(x, \frac {\pr v(x,F,q)}{\pr x}\big) \,\,
\frac\pr{\pr x_j}\,\,v (x,F,q)\,\,dx\,,
\tag{1.16}
$$
where $v(x,F,q)$ is the maximal solution of the problem (1.11), (1.12),
$C_A (F,0) = 0$. For the main properties of this capacity, in
particular the continuity with respect to $\,\,q\,$,
we refer to [5].

We assume that the following condition is satisfied.
\medskip
{\bf Condition C.} {\it There exists a function $c (x,q)$, continuous in
$x,q \in \Om \ti R^1$, such that for an arbitrary point $x \in \Om$
and an arbitrary $q \in R^1$ we have
$$
\underset{r\to 0}\to\lim \,\,\bigg\{\underset{s\to \infty}\to\lim\,\,
\frac 1{\text{\rm meas}\, K(x,r)} C_A \big( K(x,r)\setminus \Om_s,q\big)\bigg\} =
c(x,q)\,,
\tag{1.17}
$$
and the convergences to the limits in (1.17) are uniform with
respect to $\,\,q\,\,$ on any bounded interval and with respect to
$x \in \Om$.}
\medskip

The main result of the paper, proved in Section 5, is the following
theorem.
\medskip
{\bf Theorem 1.2.} {\it Assume that conditions  A.1, A.2, B, C and (1.5) are
satisfied. Let $u_s (x)$ be a sequence of solutions of the problem (0.1),
(0.2) which converges weakly in $W_m^1 (\Om)$ to a function $u_0 (x)$.
Then the function $u_0 (x)$ is a solution of the problem
$$
\underset{j=1}\to{\overset n\to\sum}\,\,\frac \pr{\pr x_j}\,\,
a_j \big(x, \frac{\pr u}{\pr x} \big) + c(x, f(x) - u(x)) =
\underset{j=1}\to{\overset n\to\sum}\,\,\frac \pr{\pr x_j}\,f_j(x),\,\,
x \in \Om,
\tag{1.18}
$$
\vskip-15pt
$$
u (x) = f (x), \,\,x \in \pr \Om,
\tag{1.19}
$$
where $c(x,q)$ is the function defined by (1.17).}
\medskip

{\bf Remark 1.3.} It is possible to establish  all results of this paper
if one replaces inequalities (1.1), (1.3) by the inequalities
$$
\underset{j=1}\to{\overset n\to\sum}\,\,a_j (x,p) p_j \> \nu_1 |p|^m,
\tag{1.20}
$$
\vskip-20pt
$$
\underset{j=1}\to{\overset n\to\sum}\,\,|a_j (x,p) -
a_j (x,q)| \< \nu_2 (|p|+|q|)^{m-2} \cdot |p-q|.
\tag{1.21}
$$
\medskip
{\bf Remark 1.4.} If we assume that conditions A.1 and B hold, 
that inequalities (1.2),
(1.20), (1.21) are satisfied  for all $x \in \Om, \,\,p,q\in R^n$, 
and that
$a_j (x,p)$ are odd and $({m-1})$-homogeneous with respect to $p$,
then the capacity defined by (1.16) satisfies the following equality
$$
C_A (F, \lm q) = |\lm|^{m-2} \lm C_A(F,q)
\tag{1.22}
$$
for every $q, \lm \in R^1$.

Under these assumptions we can formulate condition C in the following
weak form.

{\bf Condition C$^\pri$.} {\it There exists a measurable function $c(x)$
such that for almost every $x\in\Om$
$$
\aligned
&\underset{r\to 0}\to\lim \,\,
\bigg\{\underset{s\to\infty}\to{\lim \inf}\,\,
\frac 1{\text{\rm meas}\, K(x,r)} C_A \big( K(x,r)\setminus \Om_s,1\big)\bigg\} =\\
= &\underset{r\to 0}\to\lim \,\,
\bigg\{\underset{s\to\infty}\to{\lim \text{sup}}\,\,
\frac 1{\text{\rm meas}\, K(x,r)} C_A \big( K(x,r)\setminus \Om_s,1\big)\bigg\} =
c(x).
\endaligned
\tag{1.23}
$$
}

If all assumptions of this Remark are satisfied, it is still possible to
prove the result of the Theorem 1.2. For the changes in the proof we refer
to the discussions of Section 6 in [6].

\head {2. Estimates for potentials and averaging functions}
\endhead

In this section we establish some integral and pointwise estimates for the
potential functions $v (x, F, q)$ introduced in Section 1 as solutions
of problems (1.11), (1.12).

Throughout the paper we shall use the notation $C_j$, $j = 1,2,...\,$, to
indicate a constant which depends only on $n, m, \nu_1, \nu_2, A, R$,
meas $\Om$ (see (1.1), (1.3), (1.7), (1.10)).

Let us fix a compact set $\,F\,$ contained in $\Om$ and let
$v(x,q) = v(x,F,q)$. For $\mu>0$ we define the set
$$
E(\mu) = \{x\in\Om_0 : |v(x,q)| \< \mu\}.
\tag{2.1}
$$
\medskip
{\bf Lemma 2.1.} {\it Assume that conditions A.1, A.2 are satisfied
and that diam $(F) \< r$. Then there exists a constant $K_1$,
depending only on $\nu_1, \nu_2, n, m$, such that the estimate
$$
\underset {E(\mu)}\to\int \,
\bigg(1 + \bigg|\frac{\pr v (x,q)}{\pr x}\bigg|\bigg)^{m-2} \,\cdot\,
\bigg|\frac{\pr v (x,q)}{\pr x}\bigg|^2\,\,dx \<\,\,
K_1 \mu |q| (|q| + r)^{m-2}\,\,C_m (F)
\tag{2.2}
$$
holds for every $q \in R^1$ and for every $\,\,\mu>0\,$.}
\medskip

{\bf Proof.} See [6], Lemma 2.1.

\medskip

It is easy to see that the inequality
$0 \< \frac 1q v(x,q) \< 1$ holds for every $q \ne 0 $ and a.e.
$x\in \Om_0$. So we obtain an estimate of the norm of the function
$v(x,q)\,$ in $\,W_m^1 (\Om_0)$ if we put $\mu = |q|$ in (2.2).

\medskip

{\bf Theorem 2.2.} {\it Assume that conditions A.1, A.2 are satisfied,
and that $F$ is contained in a cube $K(x_0,r)$. Then there exists
a constant $K_2$, depending only on $\nu_1, \nu_2, n, m$, such that
for every $x\in {K(x_0,3r)\setminus K(x_0,r)}$ we have
$$
|v (x,q)| \< K_2 |q| \cdot \bigg[\frac r{\rho(x, K(x_0,r))}\bigg]^{n-1}
\,\,\cdot\,\,\bigg[\frac{C_m (F)}{r^{n-m}}\bigg]^{\frac 1{m-1}},
\tag{2.3}
$$
where $\rho (x,K(x_0,r))$ is the
distance from the point $\,\,x\,\,$ to the cube $K(x_0,r)$.}

\medskip

{\bf Proof.} See [10], Theorem 2.5.

\medskip

{\bf Lemma 2.3.} {\it Assume that the conditions of Theorem 2.2 and
the inequalities
$$
C_m(F) \< A\, r^n, \quad |q|^{m-1} r \< 1
\tag{2.4}
$$
are satisfied. Then there exists a constant $K_3$, depending only on
$\nu_1$, $\nu_2$, $n$, $m$ and $A$, such that the estimate
$$
|v(x,q)| \< K_3 |q| [|q| + r]^{m-2} \cdot r^2
\tag{2.5}
$$
holds for $x\in K(x_0, 2r)\setminus K(x_0, \frac{3r}2)$.}
\medskip

{\bf Proof.} We consider the case $q>0$. For $\frac r2 <\rho <r$ we
define two numerical sequences
$$
\rho_j^{(1)} = \frac{\rho}2 [1+2^{-j}],\,\,
\rho_j^{(2)} = \frac{\rho}2 [3-2^{-j}],\,\, j = 1,2,...\,,
$$
and smooth functions $\varphi_j(x)$, equal to one on the set
$G_j = K(x_0,r + \rho_j^{(2)}) \setminus  K(x_0,r + \rho_j^{(1)})$,
vanishing outside $G_{j+1}$, and such that $0 \< \varphi_j (x) \< 1,\,\,
\big|\frac{\pr \varphi_j(x)}{\pr x}\big| \< \frac 2{\rho}^{j+3}$.

Let us use the test function
$
[v(x,q)]^{\sigma +1} [\varphi_j(x)]^{\tau +m}
$
in the integral identity corresponding to the
boundary value problem (1.11), (1.12),
where $\sigma, \tau$ are arbitrary numbers greater than one.
Estimating by means of condition A.2 and Young's inequality
we obtain
$$
\aligned
&\underset {G_{j+1}}\to\int \big[1 + \big|\frac{\pr v}{\pr x}\big|\big]^{m-2}
\,\, \big|\frac{\pr v}{\pr x}\big|^2 v^\sigma\,\varphi_j^{\tau+m}\,\,dx \<\\
\< C_1 \tau^m \,&\underset {G_{j+1}}\to\int\,
\bigg[v^{\sigma +2}\big(\frac{2^j}\rho\big)^2 \varphi_j^{\tau+m-2} +
v^{\sigma+m} \big(\frac{2^j}\rho\big)^m \varphi_j^\tau\bigg] \,\,dx\,.
\endaligned
\tag{2.6}
$$

We can estimate $v(x)$ on the set $G_{j+1}$ by using
inequalities (2.3), (2.4) and we obtain
$v(x,q) \< C_2\rho$, which, together with (2.6), yields
$$
\underset {G_{j+1}}\to\int
\big|\frac{\pr v}{\pr x}\big|^2 v^\sigma\,\varphi_j^{\tau+m}\,\,dx \<
C_3 \tau^m \frac{2^{jm}}{r^2} \underset {G_{j+1}}\to\int
v^{\sigma+2} \cdot \varphi_j^\tau \,dx\,.
\tag{2.7}
$$

Define
$$
m_j = \text {ess sup} \{v(x,q) : x\in G_j\}\,.
\tag{2.8}
$$
{}From inequality (2.7) and Lemma 2.7 of [10] we obtain the following
estimate
$$
m_j^2 \< C_4 \frac{2^{\frac{jmn}2}}{r^n} \,\underset{G_{j+1}}\to\int
v^2 \varphi_j^2 dx.
\tag{2.9}
$$

The integral in the right-hand side of the last inequality is estimated
using Poincar\'e's inequality (see, e.g., [8], Chapter 8, Lemma 1.4)
and (2.2):
$$
\aligned
&\underset{G_{j+1}}\to\int \, v^2 \varphi_j^2 dx \<
\underset{G_{j+1}}\to\int \,|\min (v (x,q), m_{j+1})|^2 dx \<\\
\< C_5 r^2 & \underset {E(m_{j+1})}\to\int
\bigg|\frac{\pr v(x,q)}{\pr x}\bigg|^2 \,dx \,\<
C_6 m_{j+1} q[q+r]^{m-2}\,\cdot\, r^{2+n}.
\endaligned
\tag{2.10}
$$

By virtue of inequalities (2.9), (2.10) we have the estimate
$$
m_j^2 \< C_7\,2^{\frac{jmn}2}\,m_{j+1} q[q+r]^{m-2}\cdot r^2\,\,\,
\text {for}\,\,\,j = 1,2,...\,,
\tag{2.11}
$$
whereby, using Lemma 2.9 of [10], it follows that
$$
m_1 \< C_8 \,q\, [q+r]^{m-2} \cdot r^2.
\tag{2.12}
$$
In conclusion, we obtain estimate (2.5) from (2.8), (2.12) and the
definition of $G_1$. This completes the proof of lemma.
\medskip
We shall now state some properties of the averaging function
$u_h(x)$ defined by
$$
u_h(x) = \frac 1{h^n} \underset {R^n}\to\int K\bigg(\frac{|x-y|}h\bigg)
\,u(y) \,dy\,,
\tag{2.13}
$$
where $K(t)$ is an infinitely differentiable function, equal to zero for
$|t|\>1$, such that
$$
\underset{R^n}\to\int\,\, K (|x|)\, dx = 1
$$
and $0\< K(t) \< c(n)$ for a suitable constant $c(n)$ depending only on
$n$.

For a given positive number $\,h\,$, let us consider the family of
points $x_\al = 2 h \al\,\,$ in $\,R^n$, where
$\al = (\al_1,...,\al_n)$ is a multi-index with integer coordinates.
Let $I(h)$ be the set of multi-indices $\,\,\al\,\,$ such that
$K(x_\al, 2h) \sb \Om$ and, for every integrable function $u(x)$, let
$$
u(\al, h) = \frac 1{[2h]^n} \underset {K(x_\al,h)}\to\int
\,\,u_h (x) \,dx
$$
be the mean value of $u_h (x)$ in the cube
$K(x_\al, h)$, where $u_h (x)$ is defined by (2.13).
\medskip
{\bf Lemma 2.4.} {\it Let $\theta$ be a constant with $1\< \theta \<2$ and
let $u(x), g(x)$ be functions from the spaces $W_m^1 (\Om), L_m(\Om)$
respectively. Assume that, for some positive constant  $Q$, the
inequalities
$$
\underset{K(x_{\al},\theta h)}\to\int
|g(x)|^m \,dx \< Q h^n,\,\, \al \in I(h)
\tag{2.14}
$$
are satisfied. Then there exists a constant $K_4$, depending only on
$\,n,m\,$, such that
$$
\aligned
\underset{\al \in I(h)}\to\sum\,\,& \underset{K(x_\al,\theta h)}\to\int \,\,
|u_h (x) - u(\al,\theta h)|^m \cdot |g(x)|^m\,dx \< \\
&\< K_4 Q \cdot h^m \underset \Om\to\int
\big|\frac{\pr u(x)}{\pr x}\big|^m\,\,dx\,.
\endaligned
\tag{2.15}
$$
}
\medskip

{\bf Proof.} See [6], Lemma 2.7.
\medskip

\head {3. Proof of the Convergence Theorem}
\endhead

Let us define the sequences $\rho_s$, $\mu_s$, $\lm_s$, $\,\,s = 1,2,...\,$,
by
$$
\underset{s\to\infty}\to\lim \,\rho_s = 0,\qquad \rho_s \> r_s, \qquad
\mu_s = \big[ \ln \frac 1{\rho_s}\big]^{-1},\qquad
\lm_s = \big\{E\big(\ln \frac 1{\rho_s}\big)\big\}^{2m},
\tag{3.1}
$$
where $r_s$ is the number which appears in the condition B and
$E\big(\ln \frac 1{\rho_s}\big)$ denotes the integer part of the
number $\ln \frac 1{\rho_s}$.

We consider the subdivision of the domain $\Om$ introduced in (1.13)
and we denote
$$
K_s (\al) = K(x_\al^{(s)}, \lm_s\rho_s),\,\,\quad
K_s^\pri (\al) = K(x_\al^{(s)}, (\lm_s-2) \rho_s).
\tag{3.2}
$$
Let $q_s(x)$ be an arbitrary sequence in $L_m (\Om)$ that converges
strongly in $L_m(\Om)$ to some function $q_0(x)$.
We introduce the sets $I_s^\pri$, $I_s^{\pri\pri}$ of multi-indices by
$$
I_s^\pri = \{\al \in I_s : |q_\al^{(s)}| > 2 \mu_s\},\,\,
I_s^{\pri\pri} = \{\al \in I_s : |q_\al^{(s)}| \< 2 \mu_s\}\,,
\tag{3.3}
$$
where $q_\al^{(s)}$ is the mean value of the function $q_s(x)$
in the cube $K_s(\al)$. Let us define the functions
$w_\al^{(s)} (x)$, $\al \in I_s$, by
$$
w_\al^{(s)} (x) = v_\al^{(s)} (x, \tl q_\al^{(s)}),
\tag{3.4}
$$
where
$$
\tl q_\al^{(s)} = q_\al^{(s)}\quad\text {for} \quad \al \in I_s^\pri,\quad
\tl q_\al^{(s)} = 2 \mu_s\quad\text {for} \quad \al \in I_s^{\pri\pri}.
\tag{3.5}
$$

For an arbitrary function $g(x)$ we denote its positive part 
by $[g(x)]_+ = \max \{g(x),0\}$. We define the cut-off functions $\varphi_\al^{(s)} (x)$
by
$$
\varphi_\al^{(s)} (x) = \frac 2{\mu_\al^{(s)}}
\min\bigg\{\bigg[|w_\al^{(s)}(x)| - \frac{\mu_\al^{(s)}}2\bigg]_+,\,\,\,
\frac{\mu_\al^{(s)}}2\bigg\}\,,
\tag{3.6}
$$
where
$$
\mu_\al^{(s)} = \mu_s \cdot \max \{1, |q_\al^{(s)}|\}.
\tag{3.7}
$$
Let $G_\al^{(s)}$ be the support of the function $\varphi_\al^{(s)} (x)$.
\medskip
{\bf Lemma 3.1.} {\it Assume that conditions A.1, A.2, and B are satisfied.
Then there exists a number $s_1$ such that the inclusions
$$
G_\al^{(s)} \sb K(x_\al^{(s)},\,\,(\lm_s - 1) \rho_s)\quad
\text{for} \quad\al \in I_s
\tag{3.8}
$$
hold for $s \> s_1$.}
\medskip

The proof is analogous with the proof of Lemma 4.1 in [10].
\medskip
{\bf Lemma 3.2.} {\it Assume that conditions A.1, A.2, and B are satisfied.
Then the inequalities
$$
\text{\rm meas}\,G_\al^{(s)} \< K_5 [\lm_s \rho_s]^{m+n} \cdot
\mu_s^{1-m} \quad \text{for} \quad \al\in I_s
\tag{3.9}
$$
hold with a constant $K_5$ depending only on $\nu_1, \nu_2, n, m, A$.}
\medskip
The proof is analogous with the proof of Lemma 4.2 in [10].
\medskip
{\bf Lemma 3.3.} {\it Assume that conditions A.1, A.2, and B are satisfied, and
let $q_s(x)$ be an arbitrary sequence converging strongly in 
$L_m(\Om)$ as $\,\,s\to\infty$. Then the sequence $r_s(x)$ defined by (1.14) 
converges to zero weakly in $W_m^1(\Om)$ and
strongly in $W_p^1(\Om)$ for any $\,p<m$.}
\medskip

{\bf Proof.} We can assume that $s\>s_1\,$, where $\,s_1\,$ is defined in
Lemma 3.1. Then from inclusions (3.8) we have
$$
G_\al^{(s)} \cap G_\beta^{(s)} = \emptyset \quad\text{for}\quad
\al\ne\beta,\,\,\al, \beta \in I_s.
\tag{3.10}
$$

Let us estimate the norm of the gradient of $r_s(x)\,\,$ in
$\,\,L_m(\Om)\,\,$ for $\,\,s\,\,$ large enough such that
$$
1 \> 2\mu_s > \mu_s \>\lm_s \rho_s.
\tag{3.11}
$$
We have
$$
\aligned
\bigg\Vert\frac{\pr r_s(x)}{\pr x}\bigg\Vert^m_{L_m(\Om)} \<
C_9 \cdot \underset{\al\in I_s}\to\sum\,\,
\underset{G_\al^{(s)}}\to\int\,\,
\bigg|\frac{\pr v_\al^{(s)}(x,q_\al^{(s)})}{\pr x}\bigg|^m\,\,dx \,\,+&\\
+\,\, C_9 \cdot \underset{\al\in I_s}\to\sum\,\,\big[\mu_\al^{(s)}\big]^{-m}
\underset{\tl E_\al^{(s)}}\to\int\,\,|v_\al^{(s)} (x,q_\al^{(s)})|^m \cdot
\bigg|\frac{\pr v_\al^{(s)}(x,\tl q_\al^{(s)})}{\pr x}\bigg|^m\,\,& 
dx\,,
\endaligned
\tag{3.12}
$$
where
$
\tl E_\al^{(s)} = \big\{x\in \Om_0 : \mu_\al^{(s)}/2 \<
|v_\al^{(s)}(x,\tl q_\al^{(s)})| \< \mu_\al^{(s)}\big\}$.

The first term in the right-hand side of (3.12) is estimated by using
inequality (2.2) and condition B:
$$
\underset{\al\in I_s}\to\sum\,\,\underset{G_\al^{(s)}}\to\int\,\,
\bigg|\frac{\pr v_\al^{(s)}(x,q_\al^{(s)})}{\pr x}\bigg|^m\,\, dx\,\,
\< \,C_{10} \underset{\al\in I_s}\to\sum\,\,
\big(|q_\al^{(s)}|^m + 1\big) [\lm_s \rho_s]^n.
\tag{3.13}
$$

{}From H\"older's inequality we have
$$
|q_\al^{(s)}| = \frac 1{[2\lm_s\rho_s]^n}\,\,
\bigg|\underset{K_s(\al)}\to\int\,\,q_s (x) dx\bigg| \,\<\,
\frac 1{[2\lm_s\rho_s]^{\frac nm}}\,
\bigg\{\underset{K_s(\al)}\to\int\,\,|q_s (x)|^m\,dx\bigg\}^{\frac 1m}
\tag{3.14}
$$
and we estimate the sum in the right-hand side of (3.13) by
$$
\underset{\al\in I_s}\to\sum\,\,|q_\al^{(s)}|^m\,\,[2\lm_s\rho_s]^n\,\,\<
\,\,\underset\Om\to\int\,|q_s(x)|^m\,dx,\quad\qquad
\underset{\al\in I_s}\to\sum\,\,[2\lm_s\rho_s]^n \,\<\,\,
\text{\rm meas}\,\Om\,.
$$

Recalling the inequality
$$
|v_\al^{(s)} (x,q_\al^{(s)})| \,\< \,\mu_\al^{(s)}\quad\text{for}
\quad x\in\tl E_\al^{(s)},\,\,\al\in I_s\,,
\tag{3.15}
$$
we can estimate the second sum in the right-hand side of (3.12)
as in (3.13) and we obtain
$$
\underset\Om\to\int\,\,\bigg|\frac{\pr r_s(x)}{\pr x}\bigg|^m\,dx\,\,\<
\,\,C_{11}\,\underset\Om\to\int \big(|q_s(x)|^m + 1\big)\,dx.
\tag{3.16}
$$

Since the function $r_s(x)$ vanishes outside
$\underset{\al\in I_s}\to\cup\,\,G_\al^{(s)}$, applying H\"older's
inequality we deduce that, for $\,\,1<p<m\,$,
$$
\big\Vert\frac{\pr r_s(x)}{\pr x}\big\Vert_{L_p(\Om)}\,\,\<\,\,
\big\Vert\frac{\pr r_s(x)}{\pr x} \big\Vert_{L_m(\Om)}\,\cdot \,
\big\{\underset{\al\in I_s}\to\sum\,\text{\rm meas}\,
G_\al^{(s)}\big\}^{\frac 1p - \frac 1m}.
$$
The right-hand side of this inequality tends to zero by
(3.1), (3.9), and (3.16). 

Since, by (3.8), $r_s(x)$ has compact support in $\,\Omega\,$ 
for $\, s\> s_1\,$, the conclusions of the lemma follow from 
Poincar\'e's inequality and Rellich's compactness theorem.
\medskip
Let $\zeta_s$ be an arbitrary sequence in $R^1$ such that
$$
\underset{s\to\infty}\to\lim\,\zeta_s = 0.
\tag{3.17}
$$
Let us define the sets $I^\pri_{1,s}$, $I^\pri_{2,s}$ of
multi-indices by 
$$
I^\pri_{1,s}=\big\{\al\in I^\pri_s:\zeta_s |q_\al^{(s)}| ^{m-1}\, \< 1\big\},
\quad
I^\pri_{2,s}=\big\{\al\in I^\pri_s:\zeta_s |q_\al^{(s)}| ^{m-1}\, > 1\big\},
\tag{3.18}
$$
and denote
$$
r_{i,s}^\pri (x) \, = \,\underset{\al\in I_{i,s}^\pri}\to\sum\,\,
v_\al^{(s)} (x, q_\al^{(s)}) \,\,\varphi_\al^{(s)} (x),\quad i=1,2.
\tag{3.19}
$$
\medskip
{\bf Lemma 3.4.} {\it Assume that the conditions of Lemma 3.3 are satisfied 
and let $\zeta_s$ be an arbitrary sequence in $R^1$ satisfying (3.17).
Then the sequence $r_{2,s}^\pri (x)$ defined by (3.19) converges strongly
to zero in $W_m^1 (\Om)$.}
\medskip
{\bf Proof.} 
Define
$$
Q_s = \bigcup_{\al\in I_{2,s}^\pri}  K_s (\al)\,.
\tag{3.20}
$$
{}From (3.14)  and from
$\displaystyle
\,\,\zeta_s^{- \frac m{m-1}} \text {\rm meas}\, Q_s\,\<\,C_{12}\,
\underset{\al\in I_{2,s}^\pri}\to\sum\,\,|q_\al^{(s)}|^m\,\,[\lm_s\rho_s]^n
\,\,$
we get
$$
\text{\rm meas}\,Q_s\,\,\<\,\,C_{12}\,\,\zeta_s^{\frac m{m-1}}\,\,
\underset\Om\to\int\,\,|q_s(x)|^m\,\,dx.
\tag{3.21}
$$

As in the proof of inequality (3.16), we obtain
$$
\underset\Om\to\int \bigg|\frac{\pr r^\pri_{2,s}(x)}{\pr x}\bigg|^m\,dx
\,\<\,C_{13}\, \underset{Q_s}\to\int\,\big(|q_s(x)|^m + 
1\big)\,\,dx\,,
$$
and the convergence to zero of the right-hand side of the last inequality 
follows from (3.17), (3.21), and the assumption on the sequence $q_s(x)$.
The proof of the lemma is complete.
\medskip
{\bf Lemma 3.5.} {\it Assume that the conditions of Lemma 3.3 are satisfied.
Then the sequence
$$
r_s^{\pri\pri} (x) = \underset{\al\in I_s^{\pri\pri}}\to\sum\,\,
v_\al^{(s)} (x, q_\al^{(s)})\, \varphi_\al^{(s)} (x)
\tag{3.22}
$$
converges strongly to zero in $W_m^1 (\Om)$.}
\medskip
The proof follows immediately from the estimate
$$
\underset\Om\to\int \bigg|\frac{\pr r^{\pri\pri}_s (x)}{\pr x}\bigg|^m\,dx
\,\,\<\,\,C_{14}\,\underset{\al\in I_s^{\pri\pri}}\to\sum\,
\big(\mu_s^m + [\lm_s\rho_s]^m \big)\, [\lm_s \rho_s]^n
\,\,\<\,\,C_{14} \big(\mu_s^m + [\lm_s\rho_s]^m \big)\, 
\text{\rm meas}\,\Omega\,,
$$
that is obtained as in (3.13), using the definition of the
set $I_s^{\pri\pri}$ in (3.3).
\medskip

{\bf Proof of Theorem 1.1.} Define the sequence $\zeta_s$ by
$$
\zeta_s = \max \big\{||z_s (x)||_{L_m(\Om)},\,\,\lm_s\rho_s\big\}\,,
\tag{3.23}
$$
where $z_s (x)$ is the sequence introduced in the statement of
Theorem 1.1. Then $\zeta_s$ tends to zero as $s\to \infty$.
Let $r_{1,s}^\pri (x)$, $r_{2,s}^\pri (x)$ be the sequences defined by
(3.19) for this choice of $\zeta_s$.

Using condition A.2, Lemmas 3.3--3.5, and the assumptions on $z_s(x)$
we obtain
$$
\underset{s\to\infty}\to\lim\,\underset{j=1}\to{\overset n\to\sum}\,
\underset\Om\to\int\,\bigg[a_j\big(x, \frac{\pr r_s (x)}{\pr x}\big) -
a_j\big(x, \frac{\pr r_{1,s}^\pri (x)}{\pr x}\big) \bigg]\,\,
\frac {\pr z_s (x)}{\pr x_j}\,\,dx = 0\,,
\tag{3.24}
$$
and it is sufficient to study the behaviour of the term
$$
J_s = \underset{j=1}\to{\overset n\to\sum}\,
\underset\Om\to\int\, a_j\big(x, \frac{\pr r_{1,s}^\pri (x)}{\pr x}\big)
\,\,\frac {\pr z_s (x)}{\pr x_j}\,\,dx .
\tag{3.25}
$$

Let $\eta_\al^{(s)} (x)$ be a function of class $C_0^\infty (\Om_0)$, 
which is equal to one on $K\bigg(x_\al^{(s)},\,\frac{3 \lm_s\rho_s}2\bigg)$,
to zero outside $K (x_\al^{(s)},\,2 \lm_s\rho_s)$, and such that
$\bigg|\frac{\pr \eta_\al^{(s)} (x)}{\pr x}\bigg| \< \frac 4{\lm_s \rho_s}$.
We rewrite $J_s$ in the form
$$
J_s = \underset{j=1}\to{\overset 3\to\sum}\, J_s^{(i)}\,,
\tag{3.26}
$$
where
$$
\aligned
& J_s^{(1)} = \underset{\al\in I_{1,s}^\pri}\to\sum\,
\underset{j=1}\to{\overset n\to\sum}\,\underset{\tl K_s (\al)}\to\int\,
\bigg[a_j\big(x, \frac \pr{\pr x} (v_\al^{(s)} \varphi_\al^{(s)})\big) -
a_j\big(x, \frac {\pr v_\al^{(s)}}{\pr x}\big)\bigg]\,\,
\frac {\pr z_s (x)}{\pr x_j}\,\,dx ,\\
& J_s^{(2)} = \underset{\al\in I_{1,s}^\pri}\to\sum\,
\underset{j=1}\to{\overset n\to\sum}\,\underset{\tl K_s (\al)}\to\int\,
a_j\big(x, \frac {\pr v_\al^{(s)}}{\pr x}\big)\,\frac \pr {\pr x_j}\,
\big[\eta_\al^{(s)} (x) z_s (x)\big]\,dx, \\
& J_s^{(3)} = \underset{\al\in I_{1,s}^\pri}\to\sum\,
\underset{j=1}\to{\overset n\to\sum}\,\underset{\tl K_s (\al)}\to\int\,
a_j\big(x, \frac {\pr v_\al^{(s)}}{\pr x}\big)\,\,
\frac \pr{\pr x_j} \big[(1 - \eta_\al^{(s)} (x)) z_s (x)\big]\,\,dx;
\endaligned
\tag{3.27}
$$
here $v_\al^{(s)} = v_\al^{(s)} (x, q_\al^{(s)})$ and
$\tl K_s (\al) = K (x_\al^{(s)}, 2 \lm_s\rho_s)$.

Define
$
E_\al^{(s)} (\mu) = \big\{ x\in \tl K_s (\al):|v_\al^{(s)}
(x, q_\al^{(s)})| \< \mu \big\}$.
The function $\varphi_\al^{(s)} (x)$ is equal to one if
$|v_\al^{(s)}(x, q_\al^{(s)})| \> \mu_\al^{(s)}$, $\al\in I_s^\pri$,
and using (1.3) and H\"older's inequality we obtain the estimate
$$
\aligned
& |J_s^{(1)}| \< C_{15}\,\bigg\{\underset{\al\in I^\pri_{1,s}}\to\sum\,\,
\underset{E_\al^{(s)} (\mu_\al^{(s)})}\to\int\,\,
\bigg[1 + \big|\frac \pr{\pr x} (v_\al^{(s)} \varphi_\al^{(s)})\big| +
\big|\frac {\pr v_\al^{(s)}}{\pr x}\big|\bigg]^m\,\,dx
\bigg\}^{\frac{m-2}m} \cdot \\
&\cdot\bigg\{\underset{\al\in I^\pri_{1,s}}\to\sum\,\,
\underset{E_\al^{(s)} (\mu_\al^{(s)})}\to\int\,\,
\big|\frac \pr{\pr x} \big[v_\al^{(s)} (1 - \varphi_\al^{(s)})\big]\big|^m\,
dx \bigg\}^{\frac 1m} \cdot \bigg\{\underset\Om\to\int
\big|\frac {\pr z_s (x)}{\pr x}\big|^m\,\,dx\bigg\}^{\frac 1m}.
\endaligned
\tag{3.28}
$$

The first factor in the right-hand side of the last inequality can be
estimated from above by a constant independent of $s$. This can be 
obtained as in the proof of inequality (3.16).

We assume now that $s$ is large enough so that inequality (3.11)
is satisfied. The second factor in the right-hand side of (3.28) is
estimated using inequalities (2.2), (3.14), and
condition B. We obtain
$$
\aligned
&\underset{\al\in I^\pri_{1,s}}\to\sum\,\,
\underset{E_\al^{(s)} (\mu_\al^{(s)})}\to\int\,\,
\big|\frac \pr{\pr x} \big[v_\al^{(s)} (1 - \varphi_\al^{(s)})\big]
\big|^m\,dx \< \\
\< C_{16}& \,\,\mu_s \underset{\al\in I^\pri_{1,s}}\to\sum\,\,
 |q_\al^{(s)}|^m \,[\lm_s \rho_s]^n \<
 C_{16}\,\, \mu_s \underset\Om\to\int\,\, |q_s (x)|^m \,dx\,,
\endaligned
$$
and the right-hand side of the last inequality tends to zero as $s\to\infty$.
Taking the assumption on $z_s (x)$ into account we obtain
$$
\underset{s\to\infty}\to\lim\,\,J_s^{(1)} = 0.
\tag{3.29}
$$

The equality
$$
J_s^{(2)} = 0
\tag{3.30}
$$
follows from the definition of the functions $v_\al^{(s)} (x, q_\al^{(s)})$
(see (1.11) and (1.12)) and from the properties of $\eta_\al^{(s)} (x)$
and $z_s (x)$.

In order to estimate $J_s^{(3)}$ we remark that the inequality
$$
|v_\al^{(s)} (x, q_\al^{(s)})| \< \ov\mu_\al^{(s)},\quad
\ov\mu_\al^{(s)} = C_{17} [\lm_s \rho_s]^2\,|q_\al^{(s)}|^{m-1}
\tag{3.31}
$$
holds for $\al\in I^\pri_{1,s}$ and
$x\in\tl K_s (\al)\setminus K\big(x_\al^{(s)}, \frac{3 \lm_s\rho_s}2\big)$.
We obtain this estimate using Lemma 2.3 and (3.11), taking into 
account that $|q_\al^{(s)}|^{m-1}\cdot \lm_s\rho_s \< 1$ for
$\al\in I_{1,s}^\pri$, which implies that the second condition in (2.4) is satisfied.

By condition A.2 and H\"older's inequality we obtain the  estimate for
$J_s^{(3)}$:
$$
|J_s^{(3)}| \< C_{18}\,\frac 1{\lm_s\rho_s}\,\, \bigg\{
\underset{\al\in I^\pri_{1,s}}\to\sum\,\,
\underset{E_\al^{(s)} (\ov\mu_\al^{(s)})}\to\int\,\,
\big|\frac {\pr v_\al^{(s)}}{\pr x}\big|^m\,dx\bigg\}^{\frac{m-1}m} \cdot
\bigg\{[\lm_s\rho_s]^m\,\underset\Om\to\int\,
\big|\frac{\pr z_s}{\pr x}\big|^m\,\,dx +{}
$$
\vskip-10pt
$$
{}+ \underset\Om\to\int |z_s (x)|^m\,dx\bigg\}^{\frac 1m} + C_{18}\,
\frac 1{\lm_s\rho_s}\,\bigg\{\underset{\al\in I^\pri_{1,s}}\to\sum\,\,
\underset{E_\al^{(s)} (\ov\mu_\al^{(s)})}\to\int\,\,
\big|\frac {\pr v_\al^{(s)}}{\pr x}\big|^2\,dx\bigg\}^{\frac 12} \cdot
\tag{3.32}
$$
\vskip-10pt
$$
\cdot \big\{[\lm_s\rho_s]^2\,\underset\Om\to\int\,
\big|\frac{\pr z_s}{\pr x}\big|^2\,dx + \underset\Om\to\int\,
|z_s (x)|^2\,dx\big\}^{\frac 12},
$$
where $v_\al^{(s)} = v_\al^{(s)} (x, q_\al^{(s)})$ and
$\ov\mu_\al^{(s)}$ is defined by (3.31). In the
right-hand side of (3.22) the factors containing $z_s (x)$ can be
estimated from above by $C_{19} \zeta_s$, where $\zeta_s$ is
defined by (3.23).

In order to check the equality
$$
\underset{s\to\infty}\to\lim\,\,J_s^{(3)} = 0,
\tag{3.33}
$$
it is sufficient to establish the estimate
$$
J_s^{(4)} := \underset{\al\in I^\pri_{1,s}}\to\sum\,\,
\underset{E_\al^{(s)} (\ov\mu_\al^{(s)})}\to\int\,\,
\bigg(\big|\frac {\pr v_\al^{(s)}}{\pr x}\big|^2 +
\big|\frac {\pr v_\al^{(s)}}{\pr x}\big|^m\bigg)\,dx \< C_{20}
[\lm_s\rho_s]^2 \cdot \zeta_s^{- \frac{m-2}{m-1}}.
\tag{3.34}
$$
This inequality follows from (2.2), (3.11), (3.14), (3.18), (3.31), and
condition B:
$$
\aligned
J_s^{(4)} & \< C_{21}\, [\lm_s\rho_s]^2 \,\,
\underset{\al\in I^\pri_{1,s}}\to\sum\,\,
|q_\al^{(s)}|^{2m-2}\,\, [\lm_s\rho_s]^n \<\\
& \< C_{22} \,[\lm_s\rho_s]^2 \, \zeta_s^{- \frac{m-2}{m-1}} \cdot
\underset\Om\to\int \,|q_s (x)|^m\,dx.
\endaligned
$$
This proves inequality (3.34) and concludes the proof of the
Convergence Theorem.

\head {4. Construction and properties of test functions}
\endhead

In this section we construct special functions which belong to the
space $\overset\circ\to W^1_m (\Om_s)$ and which will be used later
as test functions in the integral identity corresponding to the boundary
value problem (0.1), (0.2).

As in Section 3 we fix the sequences
$\rho_s$, $\mu_s$, $\lm_s$ introduced in (3.1), and the subdivision of
the domain $\Om$ defined by (1.13). For $s = 1,2,...\,$ and $\al\in 
I_s$ we define $I_s (\al)$ as the set  of all multi-indices 
with integer coordinates such that
$K (2\rho_s\beta, \rho_s)\subset 
K_s (\al)\setminus \overset\circ\to K_s^\pri (\al)$,
where $K_s (\al)$, $K_s^\pri (\al)$ are the cubes defined in (3.2) and
$\overset\circ\to K_s^\pri (\al)$ is the interior of the cube
$K_s^\pri (\al)$. For $\beta\in I_s (\al)$ we set 
$x_{\al\beta}^{(s)}=2\rho_s\beta$ and 
$K_s (\al,\beta) = K (x_{\al\beta}^{(s)}, \rho_s)$. Then we have the 
following decomposition:
$$
K_s (\al)\setminus \overset\circ\to K_s^\pri (\al) =
\underset{\beta\in I_s(\al)}\to\bigcup\, K_s (\al,\beta).
\tag{4.1}
$$

Let $|I_s|$, $|I_s(\al)|$ be the numbers of multi-indices belonging to
the sets $I_s$ and $I_s (\al)$ respectively. It is easy to see that
$$
|I_s| \< C(\Om)\,[\lm_s \rho_s]^{-n},\quad
|I_s (\al)| \< 2 n \lm_s^{n-1},
\tag{4.2}
$$
where the constant $C(\Om)$ depends only on the measure of $\Om$.

Let $g(x)$ be an arbitrary function of class $C_0^\infty (\Om)$.
Let us consider the sequence
$$
q_s (x) = f_s (x) - u_0^{(s)} (x) - g (x),
\tag{4.3}
$$
where
$$
\aligned
f_s (x) = \frac 1{[\lm_s\rho_s]^n}\,\,&\underset{R^n}\to\int\,\,
K\bigg(\frac{|x-y|}{\lm_s\rho_s}\bigg)\, f(y)\,dy,\\
u_0^{(s)} (x) = \frac 1{[\lm_s\rho_s]^n}\,\,&\underset{R^n}\to\int\,\,
K\bigg(\frac{|x-y|}{\lm_s\rho_s}\bigg)\, u_0(y)\,dy,
\endaligned
$$
$f(x)$ is the boundary function from (0.2), $u_0 (x)$ is the weak
limit of the sequence $u_s (x)$, solutions of the boundary value
problem (0.1), (0.2), and the kernel $K (t)$ is the same as in (2.13).

We define new cut-off functions $\tilde\varphi_\al^{(s)} (x)$ by
$$
\tilde\varphi_\al^{(s)} (x) = \frac 2{\mu_\al^{(s)}}
\min\bigg\{\bigg[v_\al^{(s)}(x,1) - \frac{\mu_\al^{(s)}}2\bigg]_+,\,\,\,
\frac{\mu_\al^{(s)}}2\bigg\}\,,
\tag{4.4}
$$
where $v_\al^{(s)}(x,1)$ and $\mu_\al^{(s)}$ are the same as in (3.4) 
and (3.7).
In accordance with [10], we can define two sequences of
nonnegative functions $\chi_{\al\beta}^{(s)} (x)$,
$\psi_{\al\beta}^{(s)} (x)$,  for  $\al\in I_s$,
$\beta\in I_s (\al)$, such that the following properties are satisfied:

1) there exists a number $s_2$ such that the inclusions
$$
\text {supp}\,\,\chi_{\al\beta}^{(s)} \sb K\big(x_{\al\beta}^{(s)},\,
\frac {3\rho_s}2\big) \quad\text{for}\quad
\al\in I_s,\,\,\beta\in I_s(\al)
\tag{4.5}
$$
holds for $s \> s_2$ where supp  $\chi_{\al\beta}^{(s)}$ is the
support of the function $\chi_{\al\beta}^{(s)} (x)$;

2) for every point $x\in R^n$ in the sequence of numbers
$
\{\chi_{\al\beta}^{(s)} (x) : \al\in I_s,\,\,\beta\in I_s (\al)\}$
no more that $2^n$ numbers are non-zero and there exists a number $K_5$
depending only on $m, n, \nu_1, \nu_2, A$ such that the inequalities
$$
\chi_{\al\beta}^{(s)} (x) \< K_5,\quad \underset{R^n}\to\int\,\,
\bigg|\frac {\pr\chi_{\al\beta}^{(s)} (x)}{\pr x}\bigg|^m\,dx \<
K_5\,\mu_s^{1-m}\cdot \rho_s^n
\tag{4.6}
$$
holds for $s = 1,2,...,\,\,\al\in I_s,\,\,\beta\in I_s (\al)$;

3) the functions $\psi_{\al\beta}^{(s)} (x)$ are defined by the
equality
$$
\psi_{\al\beta}^{(s)} (x) = \chi_{\al\beta}^{(s)} (x)
\big\{1 - \underset{\g\in I_s}\to\sum\,\,\tilde\varphi_\g^{(s)} (x)\big\},
\quad x\in R^n;
\tag{4.7}
$$

4) the following equalities hold:
$$
\underset{\al\in I_s}\to\sum\,\, \underset{\beta\in I_s (\al)}\to\sum\,\,
\chi_{\al\beta}^{(s)} (x) = 1\quad\text{for}\quad
x\in \underset{\al\in I_s}\to\bigcup\,
\underset{\beta\in I_s (\al)}\to\bigcup\,\,
\{K_s (\al,\beta)\setminus \Om_s\}\,,
\tag{4.8}
$$
\vskip-20pt
$$
\underset{\al\in I_s}\to\sum\, \tilde \varphi_\al^{(s)} (x) +
\underset{\al\in I_s}\to\sum\,\, \underset{\beta\in I_s (\al)}\to\sum\,\,
\psi_{\al\beta}^{(s)} (x) = 1\quad\text{for}\quad
x\in \underset{\al\in I_s}\to\bigcup\,\{K_s (\al)\setminus \Om_s\}\,.
\tag{4.9}
$$

We shall assume later that
$$
s \> \max \{s_1, s_2\}.
\tag{4.10}
$$
Remark that from inclusions (3.8) and (4.5) we obtain that for every
$x\in R^n$, $\al, \g\in I_s$, $\beta\in I_s (\al)$ we have
$$
\chi_{\al\beta}^{(s)} (x)\,\, \tilde\varphi_\g^{(s)} (x) = 0\,,
\quad 
\chi_{\al\beta}^{(s)} (x)\,\, \varphi_\g^{(s)} (x) = 0
\quad \text{if}\quad
\al\ne\g.
\tag{4.11}
$$

Let us introduce the sequence
$$
h_s (x) = f (x) - q_s (x) + r_s (x) +
\underset{i=1}\to{\overset 3\to\sum} r_s^{(i)} (x),
\tag{4.12}
$$
where
$$
\aligned
r_s^{(1)} (x) &= \underset{\al\in I_s}\to\sum\,\,
[q_s (x) - q_\al^{(s)}]\,\tilde\varphi_\al^{(s)} (x),\\
r_s^{(2)} (x) &= q_s (x) \,\underset{\al\in I_s}\to\sum\,\,
\underset{\beta\in I_s(\al)}\to\sum\,\,
\psi_{\al\beta}^{(s)} (x),\\
r_s^{(3)} (x) &= \underset{\al\in I_s}\to\sum\,\,
\underset{\beta\in I_s(\al)}\to\sum\,\,
[q_\al^{(s)}\,\tilde\varphi_\al^{(s)} (x) 
- v_\al^{(s)} (x, q_\al^{(s)})\,\varphi_\al^{(s)} (x)]\,
\chi_{\al\beta}^{(s)} (x),
\endaligned
\tag{4.13}
$$
and the sequences
$r_s (x)$ and $q_s (x)$ are defined by
(1.14) and (4.3). The sequences
$q_\al^{(s)}$ and $\varphi_s^{(s)} (x)$ are the same as in Section 3, with
$q_s (x)$ defined by (4.3).
\medskip
{\bf Lemma 4.1.} {\it Assume that conditions A.1 and A.2 are 
satisfied, and let $\ov g (x)$ be an arbitrary function in the space
$C_0^\infty (\Om)$. Then there exists a number $s_3 (\ov g)$,
depending on $\ov g (x)$, such that the inclusion
$$
\ov g (x) [u_s (x) - h_s (x)] \in \overset\circ\to W_m^1 (\Om_s)
\tag{4.14}
$$
holds for $ s \> \max \{s_1, s_2, s_3 (\ov g)\}$.}
\medskip
{\bf Proof.} By the definition of the functions
$v_\al^{(s)} (x, q_\al^{(s)}),\,\,\varphi_\al^{(s)} (x)$, Lemma 3.1,
and inclusion (4.8) we obtain that the function
$$
r_s^{(4)} (x) := \underset{\al\in I_s}\to\sum\,\,
[v_\al^{(s)} (x, q_\al^{(s)})\,\varphi_\al^{(s)} (x) 
- q_\al^{(s)}\,\tilde\varphi_\al^{(s)} (x)] \cdot
\big\{1 - \,\underset{\g\in I_s}\to\sum\,\,
\underset{\delta\in I_s(\g)}\to\sum\,\,\chi_{\g\delta}^{(s)} (x)\big\}
\tag{4.15}
$$
belongs to $\overset\circ\to W_m^1 (\Om_s^\pri)$,
where $\Om_s^\pri = \Om \setminus \{\underset{\al\in I_s}\to\bigcup\,\,
[K_s (\al)\setminus \Om_s]\}$.

{}From (4.9) we obtain the inclusion
$$
r_s^{(5)} (x) := q_s (x) \big\{ 1 - \underset{\al\in I_s}\to\sum\,
\tilde\varphi_\al^{(s)} (x) - \underset{\al\in I_s}\to\sum\,
\underset{\beta\in I_s(\al)}\to\sum\,\psi_{\al\beta}^{(s)} (x)\big\}
\in \overset\circ\to W_m^1 (\Om_s^\pri).
\tag{4.16}
$$
Taking (4.11)  into account we obtain
$$
u_s (x) - h_s (x) = u_s (x) - f(x) - r_s^{(4)} (x) + r_s^{(5)} (x)
\in \overset\circ\to W_m^1 (\Om_s^\pri).
$$

Inclusion (4.14) follows now from the construction of the
subdivision (1.13) of the domain $\Om$ and from the choice of the
function $\ov g (x)$. The proof of lemma is complete.
\medskip
{\bf Lemma 4.2.} {\it Assume that conditions A.1, A.2, and B are
satisfied. Then the sequences $r_s^{(i)} (x)$, $i= 1,2,3$, defined
by (4.13), converge strongly to zero in the space $W_m^1 (\Om)$ as
$s\to\infty$.}
\medskip
{\bf Proof.} Assume that $s$ is large enough so that inequalities
(3.11) and (4.10) are satisfied. Using (2.2), (3.7), (4.4), and
condition B we have the estimate
$$
\left\Vert\frac{\pr \tilde\varphi_\al^{(s)} (x)}{\pr x}
\right\Vert^m_{L_m(\Om)}\,\<
2^m \,[\mu_\al^{(s)}]^{-m}\,\underset{\overline E_\al^{(s)}}\to\int\,
\left|\frac{\pr v_\al^{(s)} (x,1)}{\pr x}\right|^m\,dx \,\< C_{23}
\,\mu_s^{1-m} \,[\lm_s\rho_s]^n,
\tag{4.17}
$$
where $\overline E_\al^{(s)}=\big\{x\in \Om_0 : 
\mu_\al^{(s)}/2 \<
v_\al^{(s)}(x,1) \< \mu_\al^{(s)}\big\}$.

Let us estimate the norm of the gradient of $r_s^{(1)} (x)$  in
 $L_m (\Om)$:
$$
\aligned
&\bigg\Vert\frac\pr{\pr x}\,r_s^{(1)} (x)\bigg\Vert^m_{L_m (\Om)} \<
C_{24}\,\underset{\al\in I_s}\to\sum\,\underset{G_\al^{(s)}}\to\int\,
\bigg|\frac{\pr q_s (x)}{\pr x}\bigg|^m\,dx +{} \\
& {}+ C_{24}\,\underset{\al\in I_s}\to\sum\,\underset\Om\to\int\,
|g(x) - g_\al^{(s)}|^m \cdot
\bigg|\frac{\pr\tilde \varphi_\al^{(s)} (x)}{\pr x}\bigg|^m\,dx +{} \\
& {}+ C_{24}\,\underset{\al\in I_s}\to\sum\,\underset\Om\to\int\,
\{|f_s(x) - f_\al^{(s)}|^m + |u_0^{(s)}(x) - u_\al^{(s)}|^m\}
\bigg|\frac{\pr\tilde \varphi_\al^{(s)} (x)}{\pr x}\bigg|^m\,dx,
\endaligned
\tag{4.18}
$$
where $f_\al^{(s)}, u_\al^{(s)}, g_\al^{(s)}$ are the mean values of the
functions  $f_s (x), u_0^{(s)} (x), g (x)$ in the cube
$K_s (\al)$. The first term in the right-hand side of (4.18) tends to
zero as $s\to\infty$ by Lemma 3.2, the strong convergence of the
sequence $q_s (x)$ in $W_m^1 (\Om)$, and the absolute continuity of
the integral. Since the function $g(x)$ is smooth,
the second term tends to zero by (4.17) and (3.1).

Using (4.17) and Lemma 2.4, the third term in the
right-hand side of (4.18) can be estimated from above by
$$
C_{25}\,\mu_s^{1-m} \,[\lm_s \rho_s]^m\,\underset\Om\to\int\,
\left[\left|\frac{\pr f(x)}{\pr x}\right|^m +
\left|\frac{\pr u_0(x)}{\pr x}\right|^m\right]\,dx\,,
$$
which vanishes as $s\to\infty$ by (3.1). This completes
the proof of the strong convergence of $r_s^{(1)} (x)$ to zero in
$W_m^1 (\Om)$.

Let $\Cal D_{\al\beta}^{(s)}$ be the support of the function
$\psi_{\al\beta}^{(s)} (x)$. Then from (4.2), (4.5) and (4.7) we have
$$
\underset{\al\in I_s}\to\sum\,\, \underset{\beta\in I_s(\al)}\to\sum\,\,
\text {\rm meas}\,\,\Cal D_{\al\beta}^{(s)} \< C_{26} \frac 1{\lm_s}.
\tag{4.19}
$$
We will use also the estimate
$$
\underset\Om\to\int\,\,\left|\frac{\pr v_\al^{(s)} (x, q)}{\pr x}
\right|^m\,\,[\chi_{\al\beta}^{(s)} (x)]^m\,dx \,\<
C_{27}\,[\mu_s^{1-m} |q|^m + 1]\,\rho_s^n,
\tag{4.20}
$$
which follows as in the proof of inequality (4.37) in [10].
{}From (4.7), (4.11) and from inequalities (4.6), (4.20) we obtain the estimate
$$
\underset\Om\to\int\,\
\left|\frac{\pr\psi_{\al\beta}^{(s)} (x)}{\pr x}\right|^m\,dx \,\<\,C_{28}
\,\mu_s^{1-2m}\,\rho_s^n.
\tag{4.21}
$$

Let us estimate the norm of $r_s^{(2)} (x)$ in $W_m^1 (\Om)$. We rewrite
$r_s^{(2)} (x)$ in the form
$$
r_s^{(2)} (x) = \underset{\al\in I_s}\to\sum\,\,
\underset{\beta\in I_s(\al)}\to\sum\,\,
\left\{ [q_s (x) - \ov q_\al^{(s)}] + \ov q_\al^{(s)} \right\}\,\,
\psi_{\al\beta}^{(s)} (x),
$$
where $\ov q_\al^{(s)}$ is the mean value of the function $q_s (x)$
in the cube $\tl K_s (\al) = K(x_\al^{(s)}, 2\lm_s\rho_s)$.
Using (4.6), (4.7), and (4.21) we obtain the inequality
$$
\aligned
\underset\Om\to\int\,
\left|\frac{\pr r_s^{(2)} (x)}{\pr x}\right|^m\,dx \,&\<\,C_{29}\,
\underset{\al\in I_s}\to\sum\,\,\underset{\beta\in I_s(\al)}\to\sum\,\,
\underset{\Cal D_{\al\beta}^{(s)}}\to\int\,
\left|\frac{\pr q_s (x)}{\pr x}\right|^m\,dx \,+{} \\
{}+ \,C_{29}\,\mu_s^{1-2m} & \underset{\al\in I_s}\to\sum\,\,
\underset{\beta\in I_s(\al)}\to\sum\,\,
\left|\ov q_\al^{(s)}\right|^m \cdot \rho_s^n +{} \\
{}+ C_{29} \,\underset{\al\in I_s}\to\sum\,
\underset{\tl K_s (\al)}\to\int &\,\left|q_s (x) - \ov q_\al^{(s)}\right|^m\,
\underset{\beta\in I_s(\al)}\to\sum\,\,
\left|\frac{\pr \psi_{\al\beta}^{(s)} (x)}{\pr x}\right|^m\,dx .
\endaligned
\tag{4.22}
$$

The first term in the right-hand side of (4.22) tends to zero as
$s\to\infty$ by (4.19) and the strong convergence of the sequence
$q_s (x)$ in $W_m^1 (\Om)$. Using (4.2) and (3.14), the
second term in the right-hand side of (4.22) is estimated from above by
$$
C_{30} \,\lm_s^{-1}\mu_s^{1-2m}\,\underset{\al\in I_s}\to\sum\,
\left|\ov q_\al^{s)}\right|^m\,[\lm_s \rho_s]^n\,\<
C_{31} \lm_s^{-1}\mu_s^{1-2m}\,\underset\Om\to\int\,
\left| q_s (x)\right|^m\,dx,
$$
which tends to zero by the choice of $\mu_s,\lm_s$.

Using Lemma 2.4 and inequalities (4.2), (4.21), 
the third term in the right-hand side of (4.22) is estimated from above  by
$$
\aligned
C_{32}\, \lm_s^{-1}\,\mu_s^{1-2m} & [\lm_s \rho_s]^m
\bigg\{\underset\Om\to\int\left[\left|\frac{\pr f(x)}{\pr x}\right|^m +
\left|\frac{\pr u_0(x)}{\pr x}\right|^m\right]\,dx +{} \\
{}+ & \underset {x\in\Om}\to\max\,\left|\frac{\pr g(x)}{\pr x}\right|^m \cdot
\,\text{\rm meas}\,\Om\bigg\},
\endaligned
$$
which converges to zero by (3.1). This concludes the proof of the 
strong convergence to zero of $r_s^{(2)} (x)$ in $W_m^1 (\Om)$.

The same property for $r_s^{(3)} (x)$
follows from the inequality
$$
\aligned
\underset\Om\to\int\,\left|\frac{\pr r_s^{(3)} (x)}{\pr x}\right|^m\,& dx
\,\<\, C_{33}\,\underset{\al\in I_s}\to\sum\,
\underset{\beta\in I_s (\al)}\to\sum\,
\left\{|q_\al^{(s)}|^m \,\mu_s^{1-2m} + 1\right\}\,\rho_s^n \,\< \\
& \< C_{34} \lm_s^{-1} \underset\Om\to\int\,
\left\{\mu_s^{1-2m}\,|q_s (x)|^m + 1\right\}\,dx\,,
\endaligned
$$
which can be obtained by using (4.2), (4.6), (4.20), and (3.14). This completes the
proof of Lemma 4.2.
\medskip
Let $g(x)$ be the same function as before, and let $g_\al^{(s)}$ be its
mean value in the cube $K_s (\al)$. We introduce the sequence
$$
g_s (x) = g(x) + \rho_s (x) +
\underset{i=1}\to{\overset 3\to\sum}\,\rho_s^{(i)} (x),
\tag{4.23}
$$
where
$$
\aligned
&\rho_s (x) = - \underset{\al\in I_s}\to\sum\,\,\frac 1{\tl q_\al^{(s)}}
\,v_\al^{(s)} (x, \tl q_\al^{(s)})\,\varphi_\al^{(s)} (x) \,g_\al^{(s)},\\
&\rho_s^{(1)}(x) = \underset{\al\in I_s}\to\sum\,\,[g_\al^{(s)} - g(x)]\,
\tilde \varphi_\al^{(s)} (x),\\
&\rho_s^{(2)}(x) = - g(x)\underset{\al\in I_s}\to\sum\,
\underset{\beta\in I_s}\to\sum\,\psi_{\al\beta}^{(s)}(x),\\
&\rho_s^{(3)}(x) = - \underset{\al\in I_s}\to\sum\,
\underset{\beta\in I_s(\al)}\to\sum\,
\left[\tilde\varphi_\al^{(s)} (x) - \frac 1{\tl q_\al^{(s)}}
 v_\al^{(s)} (x, \tl q_\al^{(s)})\, \varphi_\al^{(s)} (x) \right]\,
g_\al^{(s)}\,\chi_{\al\beta}^{(s)}(x).
\endaligned
\tag{4.24}
$$
Here $\varphi_\al^{(s)}(x)$, $\tilde\varphi_\al^{(s)} (x)$, 
$\psi_{\al\beta}^{(s)}(x)$,
$\chi_{\al\beta}^{(s)}(x)$ are the same functions as in (4.13),
$$
\tl q_\al^{(s)} = q_\al^{(s)}\quad \text{for}\quad
\al\in I_s^\pri,\quad
\tl q_\al^{(s)} = 2\mu_s \quad \text{for}\quad
\al\in I_s^{\pri\pri},
$$
and $q_\al^{(s)}$ is the mean value in the cube
$K_s (\al)$ of the function $q_s (x)$ defined by (4.3).
\medskip
{\bf Lemma 4.3.} {\it Assume that conditions of Lemma 4.1 are satisfied.
Then there exists a number $s_4 (\ov g)$ depending on $\ov g (x)$
such that
$$
\ov g (x)\,g_s (x) \in \overset\circ\to W_m^1 (\Om_s)
\tag{4.25}
$$
for $s \> s_4 (\ov g)$.}
\medskip
The proof is analogous with the proof of Lemma 4.1.
\medskip
{\bf Lemma 4.4.} {\it Assume that conditions A.1, A.2, and B are satisfied.
Then the sequence $\rho_s (x)$ is bounded in $W_m^1 (\Om)$ and
converges to zero strongly in $W_p^1 (\Om)\,\,$ for $\,\,p<m$.}
\medskip
The proof is analogous with the proof of Lemma 3.3.
\medskip
{\bf Lemma 4.5.} {\it Assume that conditions A.1, A.2, and B are satisfied.
Then the sequences $\rho_s^{(i)} (x)$, $i = 1,2,3$, converge strongly
to zero in $W_m^1 (\Om)\,\,$ as $\,\,s\to\infty$.}
\medskip
The proof is analogous with the proof of Lemma 4.2.
\medskip

\head {5. Construction of the limit boundary value problem}
\endhead

Using condition C we can conclude that for an arbitrary positive
number $\ve$ there exist two positive numbers $r (\ve)$ and $s (\ve)$, and
a sequence $r_s (\ve)$ converging to zero as $s\to\infty$, such that
the inequality
$$
\left|\frac 1{\text{\rm meas}\, K (x,r)}\, C_A (K (x,r)\setminus \Om_s, q) -
c (x,q)\right| < \ve
\tag{5.1}
$$
holds for $s\> s (\ve),\quad r_s (\ve) \< r \< r (\ve),\quad
|q| \< \frac 1 \ve, \quad x\in\Om$.
\medskip
{\bf Lemma 5.1.} {\it Assume that conditions A.1, A.2, B, and C are
satisfied. Then the function $c(x,q)$ defined by condition C
satisfies the inequality
$$
|c(x,q)| \< K_6 |q|^{m-1}
\tag{5.2}
$$
with a constant $K_6$ depending only on $n, m, \nu_1, \nu_2, A$.}
\medskip
{\bf Proof.} Inequality (5.2) follows immediately from definition
of the function $c(x,q)$ and inequality (2.2).
\medskip
{\bf Proof of Theorem 1.2.}
Let us fix a positive number $\ve$ and let $\rho_s$ be the sequence 
defined by
$$
\rho_s = r_s + r_s (\ve).
\tag{5.3}
$$
Let $\mu_s$ and $\lm_s$ be the sequences defined by equalities (3.1) 
with this particular choice of $\rho_s$.

We fix some function $g(x)$ in the space $C_0^\infty (\Om)$ and choose
a nonnegative function $\ov g (x) \in C_0^\infty (\Om)$ such that
$$
\ov g (x)\,g (x) = g (x),\qquad |g (x)| \< 1,\qquad
\ov g(x) \< 1\qquad\text{for}\quad x\in\Om.
\tag{5.4}
$$

Using inequality (1.2) we get
$$
0 \< \underset{j=1}\to{\overset n\to\sum}\,\underset{\Om_s}\to\int\,
\left[a_j\left(x, \frac{\pr u_s}{\pr x}\right) -
a_j\left(x, \frac{\pr h_s}{\pr x}\right)\right]\,
\frac{\pr(u_s - h_s)}{\pr x_j}\,\ov g (x)\,dx =
\underset{i=1}\to{\overset 4\to\sum}\, R_s^{(i)}\,,
\tag{5.5}
$$
where $u_s (x)$ is the solution of problem (0.1), (0.2),
$h_s (x)$ in the function introduced in (4.12), and $R_s^{(i)}$ are
defined by equalities
$$
\aligned
R_s^{(1)} &= -\underset{j=1}\to{\overset n\to\sum}\,\underset{\Om_s}\to\int\,
\left[a_j\left(x, \frac{\pr u_s}{\pr x}\right) -
a_j\left(x, \frac{\pr h_s}{\pr x}\right)\right]\,
(u_s - h_s + g)\,\frac {\pr\ov g (x)}{\pr x_j}\,dx,\\
R_s^{(2)} &= \underset{j=1}\to{\overset n\to\sum}\,\underset{\Om_s}\to\int\,
a_j\left(x,\frac{\pr u_s}{\pr x}\right) \frac\pr{\pr x_j}\,
\left\{\ov g (x) [u_s (x) - h_s (x)] \right\}\,dx,\\
R_s^{(3)} &= -\underset{j=1}\to{\overset n\to\sum}\,\underset{\Om_s}\to\int\,
\left[a_j\left(x,\frac{\pr h_s}{\pr x}\right) -
a_j\left(x,\frac{\pr r_s}{\pr x}\right)\right]\,\frac\pr{\pr x_j}\,
\left\{\ov g (x) [u_s (x) - h_s (x)] \right\}\,dx,\\
R_s^{(4)} &= -\underset{j=1}\to{\overset n\to\sum}\,\underset{\Om_s}\to\int\,
a_j\left(x,\frac{\pr r_s}{\pr x}\right) \frac\pr{\pr x_j}\,
\left\{\ov g (x) [u_s (x) - h_s (x)] \right\}\,dx.
\endaligned
\tag{5.6}
$$

By Lemmas 3.3 and 4.2 the sequence
$u_s (x) - h_s (x) + g (x)$ converges to zero strongly in $L_m (\Om)$ and 
then by the same Lemmas and (3.1) we obtain the equality
$$
\underset{s\to\infty}\to\lim\, R_s^{(1)} = 0.
\tag{5.7}
$$

Using the definition of $u_s (x)$ and Lemmas 3.3, 4.1, 4.2 we get the
equality
$$
\underset{s\to\infty}\to\lim\, R_s^{(2)} =
- \underset{j=1}\to{\overset n\to\sum}\,\underset{\Om}\to\int\,
f_j (x) \,\frac{\pr g (x)}{\pr x_j}\,dx.
\tag{5.8}
$$

We check now that
$$
\underset{s\to\infty}\to\lim\, R_s^{(3)} =
\underset{j=1}\to{\overset n\to\sum}\,\underset{\Om}\to\int\,
a_j\left(x,\frac{\pr u_0}{\pr x} + \frac{\pr g}{\pr x}\right)\,
\frac {\pr g}{\pr x} \,dx.
\tag{5.9}
$$
By inequality (1.3) and Lemmas 3.3, 4.2 we obtain that the sequence
$$
b_j^{(s)} (x) = a_j \left(x, \frac{\pr h_s (x)}{\pr x}\right) -
a_j \left(x, \frac{\pr r_s (x)}{\pr x}\right)
$$
converges in measure to the function $a_j\left(x,\frac{\pr u_0}{\pr x} +
\frac{\pr g}{\pr x}\right)$. To prove that this sequence
converges strongly in $L_{\frac m{m-1}} (\Om)$, it is sufficient to
establish that the integrals
$$
\underset\Om\to\int\,\left|b_j^{(s)} (x)\right|^{\frac m{m-1}}\,dx,\,\,\,
s = 1,2,...
\tag{5.10}
$$
satisfy the absolute continuity property uniformly with respect to $s$.

Using (1.3) and H\"older's inequality we obtain the estimate
$$
\aligned
\underset E\to\int \left|b_j^{(s)} (x)\right|^{\frac m{m-1}}\,& dx \,
\< C_{35}\left\{ \underset E\to\int\,
\left[1 + \left|\frac{\pr h_s (x)}{\pr x}\right| +
\left|\frac{\pr r_s (x)}{\pr x}\right|\right]^m\,dx\,
\right\}^{\frac{m-2}{m-1}} \cdot\\
&\cdot \left\{ \underset E\to\int\, \left|\frac{\pr h_s (x)}{\pr x} -
\frac{\pr r_s (x)}{\pr x}\right|^m\,dx\,\right\}^{\frac 1{m-1}}
\endaligned
$$
for an arbitrary subset $E$ of the set $\Om$. The last inequality and Lemma 4.2
guarantee the uniform absolute continuity for the sequence
of integrals (5.10), and hence the strong convergence of $b_j^{(s)} (x)$ in 
$L_{\frac m{m-1}} (\Om)$. Using this property 
and Lemmas 3.3, 4.2 we obtain equality (5.9).

We transform the term $R_s^{(4)}$ in the following way:
$$
R_s^{(4)} = \underset{i=5}\to{\overset 9\to\sum}\,R_s^{(i)}\,,
\tag{5.11}
$$
where
$$
\aligned
R_s^{(5)} &= -\underset{j=1}\to{\overset n\to\sum}\,
\underset\Om\to\int\, a_j\left(x,\frac{\pr r_s}{\pr x}\right)\,
\frac \pr{\pr x_j}\,\left\{\ov g(x) [u_s(x) - h_s(x) + g_s(x)]\right\}\,dx,\\
R_s^{(6)} &= \underset{j=1}\to{\overset n\to\sum}\,
\underset\Om\to\int\, a_j\left(x,\frac{\pr r_s}{\pr x}\right)\,
\frac {\pr g(x)}{\pr x_j}\,dx,\\
R_s^{(7)} &= \underset{j=1}\to{\overset n\to\sum}\,
\underset\Om\to\int\, a_j\left(x,\frac{\pr r_s}{\pr x}\right)\,
\frac {\pr \ov g(x)}{\pr x_j}\,\cdot\rho_s (x)\,dx,\\
R_s^{(8)} &= \underset{j=1}\to{\overset n\to\sum}\,
\underset\Om\to\int\, a_j\left(x,\frac{\pr r_s}{\pr x}\right)\,
\frac {\pr \rho_s(x)}{\pr x_j}\,\cdot \ov g(x)\,dx,\\
R_s^{(9)} &= \underset{i=1}\to{\overset 3\to\sum}\,
\underset{j=1}\to{\overset n\to\sum}\,\underset\Om\to\int\,
a_j\left(x,\frac{\pr r_s}{\pr x}\right)\,\frac \pr{\pr x_j}\,
\left[\ov g(x) \rho_s^{(i)}(x)\right]\,dx,
\endaligned
$$
and the functions $g_s(x), \rho_s(x), \rho_s^{(i)}(x)$ are defined 
by (4.23), (4.24).

By virtue of Lemmas 3.3, 4.1--4.5 the sequence
$
z_s(x) = \ov g(x) [u_s(x) - h_s(x) + g_s(x)]
$
satisfies the following conditions for $\,\,s\,\,$ large enough:
$z_s(x)\in\overset\circ\to W_m^1 (\Om_s)$ and $z_s(x)$ converges
to zero weakly in $W_m^1(\Om)$. Then by the Convergence Theorem 1.1 we obtain
$$
\underset{s\to\infty}\to\lim\, R_s^{(5)} = 0.
\tag{5.12}
$$
Using (1.3) and Lemmas 3.3, 4.4. 4.5 we have
$$
\underset{s\to\infty}\to\lim\, R_s^{(6)} =
\underset{s\to\infty}\to\lim\, R_s^{(7)} =
\underset{s\to\infty}\to\lim\, R_s^{(9)} = 0.
\tag{5.13}
$$

It remains to study  the behaviour of $R_s^{(8)}$. {}From the
definitions of $r_s(x), \rho_s(x)$ we obtain for $s \> s_1$:
$$
R_s^{(8)} = - \underset{i=10}\to{\overset 13\to\sum}\,R_s^{(i)}\,,
\tag{5.14}
$$
where
$$
\aligned
R_s^{(10)} = \underset{\al\in I_s}\to\sum\,\,
\underset{j=1}\to{\overset n\to\sum}\,&
\frac {g_\al^{(s)}}{\tl q_\al^{(s)}}\,\underset{\Om_0}\to\int\,
a_j\left(x,\frac{\pr\left(v_\al^{(s)} \varphi_\al^{(s)}\right)}{\pr x}
\right) \cdot \frac \pr{\pr x_j}
\left(\tl v_\al^{(s)} \varphi_\al^{(s)}\right) \cdot
\left[\ov g(x) - \ov g_\al^{(s)}\right]\,dx,\\
R_s^{(11)} = \underset{\al\in I_s}\to\sum\,\,
\underset{j=1}\to{\overset n\to\sum}\, &
\frac {g_\al^{(s)} \ov g_\al^{(s)}}{\tl q_\al^{(s)}}\,
\underset{\Om_0}\to\int\,
\bigg\{a_j\left(x,\frac{\pr\left(v_\al^{(s)} \varphi_\al^{(s)}\right)}{\pr x}
\right) \cdot \frac \pr{\pr x_j}
\left(\tl v_\al^{(s)} \varphi_\al^{(s)}\right) - \\
& - a_j \left(x,\frac{\pr v_\al^{(s)}}{\pr x}\right)\,
\frac{\pr \tl v_\al^{(s)}}{\pr x_j}\bigg\}\,\,dx,\\
R_s^{(12)} = \underset{\al\in I_s^\pri}\to\sum\,\,
\underset{j=1}\to{\overset n\to\sum}\, &
\frac {g_\al^{(s)} \ov g_\al^{(s)}}{q_\al^{(s)}}\,
\underset{\Om_0}\to\int\, a_j\left(x,\frac{\pr v_\al^{(s)}}{\pr x}\right)
\frac{\pr v_\al^{(s)}}{\pr x_j}\,dx,\\
R_s^{(13)} = \underset{\al\in I_s^{\pri\pri}}\to\sum\,\,
\underset{j=1}\to{\overset n\to\sum}\, &
\frac {g_\al^{(s)} \ov g_\al^{(s)}}{2\mu_s}\,
\underset{\Om_0}\to\int\, a_j\left(x,\frac{\pr v_\al^{(s)}}{\pr x}\right)
\,\frac{\pr \tl v_\al^{(s)}}{\pr x_j}\,dx;
\endaligned
$$
here
$v_\al^{(s)} = v_\al^{(s)} (x, q_\al^{(s)})$,
$\tl v_\al^{(s)} = v_\al^{(s)} (x, \tl q_\al^{(s)})$, $\ov g_\al^{(s)}$
is the mean value of the function $\ov g(x)$ in the cube
$K_s(\al)$, and the sets $I_s^\pri$, $I_s^{\pri\pri}$ are defined by (3.3).

Using inequalities (1.3), (2.2), (3.14), (3.15), condition B,
and the smoothness of the function $\ov g(x)$, we obtain
$$
\underset{s\to\infty}\to\lim\,R_s^{(10)} = 0.
\tag{5.15}
$$
We check now that
$$
\underset{s\to\infty}\to\lim\,R_s^{(11)} = 0, \quad
\underset{s\to\infty}\to\lim\,R_s^{(13)} = 0.
\tag{5.16}
$$
The first equality in (5.16) is established as in the
proof of equality (3.29). We only need to observe that, by
(2.2) and condition B, we have the estimate
$$
\underset{\al\in I_s}\to\sum\,\frac 1{|\tl q_\al^{(s)}|^m}\,
\underset{\Om_0}\to\int\,
\left\{\left|\frac{\pr v_\al^{(s)}(x,q_\al^{(s)})}{\pr x}\right|^m +
\left|\frac{\pr v_\al^{(s)}(x,\tl q_\al^{(s)})}{\pr x}\right|^m\right\}\,
dx \< C_{36}.
$$
The second equality in (5.16) follows immediately from estimate (2.2)
and from the definition of $I_s^{\pri\pri}$.

In the rest of the paper we use the notation
$\delta_i(t)$, $\g_i(t)$, $i = 1,2,...\,$, to indicate
nonnegative functions on $R^1$ satisfying the conditions
$$
\underset{t\to 0}\to\lim\, \delta_i (t) = 0, \quad
\underset{t\to\infty}\to\lim\, \g_i (t) = 0.
\tag{5.17}
$$
\medskip
{\bf Lemma 5.2.} {\it Assume that conditions A.1, A.2, B, and C are satisfied.
Then there exist functions $\delta_1(t),\, \g_1(t)$ satisfying
conditions (5.17) such that
$$\aligned
&\bigg|\underset{\al\in I_s^\pri}\to\sum\,\,
\underset{j=1}\to{\overset n\to\sum}\,
\frac{\ov g_\al^{(s)} g_\al^{(s)}}{q_\al^{(s)}}\,
\underset{\Om_0}\to\int\, a_j\left(x,\frac{\pr v_\al^{(s)}}{\pr x}\right)\,
\frac{\pr v_\al^{(s)}}{\pr x_j}\,dx -{}\\
& {}- \underset\Om\to\int\, c(x,q_0(x))\,g(x)\,dx\bigg|
\< \g_1 (s) + \delta_1 (\ve)\,,
\endaligned
\tag{5.18}
$$
where $\ve$ is the number fixed in the definition of $\rho_s$ in
(5.3), and  $q_0(x) = f(x) - u_0(x) - g(x)$.}
\medskip

{\bf Proof.} Define the sets of multi-indices
$$
I_s(\ve) = \left\{\al\in I_s^\pri : |q_\al^{(s)}| \< \frac 1\ve\right\},\quad
J_s(\ve) = \left\{\al\in I_s^\pri : |q_\al^{(s)}| > \frac 
1\ve\right\},
$$
and let $Q_s(\ve)$ be the union of all cubes $K_s(\al)$ with $\al\in J_s(\ve)$.
As in (3.16) and (3.21) we obtain the estimates
$$
\text{\rm meas}\,Q_s(\ve) \< C_{37}\,\ve^m \underset\Om\to\int\,
|q_s (x)|^m\,dx,
\tag{5.19}
$$
$$
\aligned
&\left|\underset{\al\in J_s(\ve)}\to\sum\,\,
\underset{j=1}\to{\overset n\to\sum}\,\,
\frac {\ov g_\al^{(s)} g_\al^{(s)}}{q_\al^{(s)}}\,
\underset{\Om_0}\to\int\,a_j \left(x,\frac{\pr v_\al^{(s)}}{\pr x}\right)\,
\frac{\pr v_\al^{(s)}}{\pr x_j}\,dx\right| \< \\
\< C_{38} \underset{Q_s(\ve)}\to\int\,&|q_0(x)|^{m-1}\,dx +
C_{38} \underset\Om\to\int\,|q_s(x)- q_0(x)|^{m-1}\,dx \<
\delta_2(\ve) + \g_2(s).
\endaligned
\tag{5.20}
$$

Using notation (1.16) and inequality (5.1) we have the estimate
$$
\aligned
&\bigg|\underset{\al\in I_s(\ve)}\to\sum\,\,
\underset{j=1}\to{\overset n\to\sum}\,\,
\frac {\ov g_\al^{(s)} g_\al^{(s)}}{q_\al^{(s)}}\,
\underset{\Om_0}\to\int\,a_j \left(x,\frac{\pr v_\al^{(s)}}{\pr x}\right)\,
\frac{\pr v_\al^{(s)}}{\pr x_j}\,dx - {}\\
{}- &\underset{\al\in I_s(\ve)}\to\sum\, \ov g_\al^{(s)} g_\al^{(s)}
c(x_\al^{(s)}, q_\al^{(s)})\, \text{\rm meas}\, K_s^\pri (\al)\bigg|
< \ve \,\,\text{\rm meas}\, \Om\,,
\endaligned
\tag{5.21}
$$
provided $s$ is so large that $\lm_s \rho_s < r (\ve)$.

In view of the continuity of the functions $c(x,q), g(x), \ov g(x)$ we
obtain the inequality
$$
\bigg|\underset{\al\in I_s(\ve)}\to\sum\, \ov g_\al^{(s)} g_\al^{(s)}
c(x_\al^{(s)}, q_\al^{(s)}) \,\text{\rm meas}\, K_s^\pri (\al) -
\underset{\al\in I_s(\ve)}\to\sum\,\underset{K_s^\pri(\al)}\to\int\,
c(x, q_\al^{(s)}) g(x) dx\bigg| \< \g_3(s).
\tag{5.22}
$$

Using inequalities (5.2), (3.14), and (5.19) we have the following
estimates
$$
\bigg|\underset{\al\in J_s(\ve)}\to\sum\,\underset{K_s^\pri(\al)}\to\int\,
c(x, q_\al^{(s)}) g(x) dx\bigg| \< C_{39} \ve \underset{\Om}\to\int\,
\left|q_s(x)\right|^m dx \< \delta_3(\ve),
\tag{5.23}
$$
$$
\bigg|\underset{\al\in I_s^{\pri\pri}}\to\sum\,
\underset{K_s^\pri(\al)}\to\int\,
c(x, q_\al^{(s)}) g(x) dx\bigg| \< C_{40}\,\mu_s^{m-1}\,
\text{\rm meas}\, \Om \< \g_4(s).
\tag{5.24}
$$
{}From the last two estimates we obtain the inequality
$$
\bigg|\underset{\al\in I_s(\ve)}\to\sum\,\underset{K_s^\pri(\al)}\to\int\,
c(x, q_\al^{(s)})\, g(x)\, dx - \underset\Om\to\int\,
c\left(x, q_s^\pri (x)\right)\, g(x)\,dx \bigg| \< \delta_3(\ve) + \g_4(s)\,,
\tag{5.25}
$$
where
$$
q_s^\pri(x) = \underset{\al\in I_s}\to\sum\,q_\al^{(s)}\,
\chi (K_s^\pri (\al))\,,
\tag{5.26}
$$
and $\chi (K_s^\pri (\al))$ is the characteristic function of the set
$K_s^\pri(\al)$.

We check that the sequence $q_s^\pri(x)$ defined by (5.26) converges
to $q_0(x)$ strongly in $L_m(\Om)$. Using Poincar\'e's inequality and
(4.2) we have
$$
\aligned
& \underset\Om\to\int \left|q_s^\pri(x) - q_0(x)\right|^m\,dx \<
C_{41}\,\,\bigg\{\underset{U_s}\to\int |q_0(x)|^m\,dx +{} \\
&{}+ \underset\Om\to\int |q_s(x) - q_0(x)|^m\,dx +
\underset{\al\in I_s}\to\sum\,\,
\underset{K_s(\al)\setminus K_s^\pri(\al)}\to\int\,\,|q_0(x)|^m\,dx +{} \\
& {}+ (\lm_s\rho_s)^m \,\underset\Om\to\int
\left|\frac{\pr q_s(x)}{\pr x}\right|^m\,dx\bigg\}\,\<\,\,\g_5(s).
\endaligned
\tag{5.27}
$$

{}From (5.2), (5.27) and the continuity of the function $c(x,q)$ we
obtain the estimate
$$
\bigg|\underset\Om\to\int \left[c(x,q_s^\pri (x)) -
c(x,q_0(x)\right] \,g(x)\, dx\bigg| \< \g_6(s).
\tag{5.28}
$$
Now inequality (5.18) follows from (5.20)--(5.22), (5.25), (5.28)
and the proof of the lemma is complete.
\medskip
Inequality (5.5), together with (5.7)--(5.9), (5.11)--(5.16), and (5.18),
implies
$$
\aligned
& \underset{j=1}\to{\overset n\to\sum}\,\underset\Om\to\int\,
\left[a_j \left(x,\frac{\pr u_0}{\pr x} +\frac{\pr g}{\pr x}\right) -
f_j (x)\right] \frac{\pr g(x)}{\pr x_j}\,\,dx - {}\\
&{}- \underset\Om\to\int c(x, f(x) - u_0(x) - g(x)) g(x) dx \<
\g_7(s) + \delta_4(\ve).
\endaligned
\tag{5.29}
$$
In (5.29) the left-hand side is independent of $\,s\,$
and $\,\ve\,$ while the right-hand side can be made arbitrarily small for
sufficiently large $\,s\,$ and sufficiently small $\,\ve$. Hence,
we obtain the inequality
$$
\aligned
&\underset{j=1}\to{\overset n\to\sum}\,\underset\Om\to\int\,
\left[a_j \left(x,\frac{\pr u_0}{\pr x} +\frac{\pr g}{\pr x}\right) -
f_j (x)\right] \frac{\pr g(x)}{\pr x_j}\,\,dx -{}\\
&{}- \underset\Om\to\int c(x, f(x) - u_0(x) - g(x)) g(x) dx \< 0
\endaligned
\tag{5.30}
$$
for an arbitrary function $g(x)\in C_0^\infty (\Om)$.

In (5.30) we can replace $g(x)$ by $\lm g(x)$, with $\lm > 0$.
Dividing both terms of (5.30) by $\lm$ and passing to the limit as
$\lm\to 0$ we obtain
$$
\underset{j=1}\to{\overset n\to\sum}\,\underset\Om\to\int\,
\left[a_j \left(x,\frac{\pr u_0}{\pr x}\right) -
f_j (x)\right] \frac{\pr g(x)}{\pr x_j}\,dx -
\underset\Om\to\int \,c(x, f(x) - u_0(x)) g(x) dx \> 0.
\tag{5.31}
$$
This inequality is true for both functions $g(x)$ and $-g(x)$ and
consequently the left-hand side of (5.31) is equal zero for
$g(x)\in C_0^\infty (\Om)$. By an approximation argument we obtain the equality
$$
\underset{j=1}\to{\overset n\to\sum}\,\underset\Om\to\int\,
\left[a_j \left(x,\frac{\pr u_0}{\pr x}\right) -
f_j (x)\right] \frac{\pr g(x)}{\pr x_j}\,dx -
\underset\Om\to\int \,c(x, f(x) - u_0(x)) g(x) dx = 0
$$
for an arbitrary function $g(x)\in\overset\circ\to W_m^1 (\Om)$. Thus 
we have established that $u_0(x)$ is a solution of equation (1.18). The
inclusion $u_0(x)\in f(x) + \overset\circ\to W_m^1 (\Om)$ follows
immediately from $u_s(x)\in f(x) + \overset\circ\to W_m^1 (\Om)$.
This completes the proof of Theorem 1.2.

\head {References}\endhead
\frenchspacing

\smallskip
\item{[1]}Casado Diaz J., Garroni A.: Asymptotic behaviour of nonlinear 
elliptic systems on varying domains. {\it SIAM J. Math. Anal.\/} {\bf 31} 
(2000), 581-624.
\smallskip
\item{[2]}Dal Maso G., Defranceschi A.: Limits of nonlinear Dirichlet problems
in varying domains. {\it Manuscripta Math\/}. {\bf 61} (1988), 251-278.
\smallskip
\item{[3]}Dal Maso G., Murat F.: Dirichlet problems in perforated domaims for
homogeneous monotone operators on $H_0^1$. {\it Calculus of Variations,
Homogenization and Continuum Mechanics (CIRM - Luminy, Marseille, 
1993)\/}, 177-202, {\it World Scientific, Singapore\/}, 1994.
\smallskip
\item{[4]}Dal Maso G., Murat F.: Asymptotic behaviour and correctors for
Dirichlet problems in perforated domains with homogeneous monotone
operators. {\it Ann. Scuola Norm. Sup. Pisa Cl. Sci. (4)\/} {\bf 24} (1997),
239-290.
\smallskip
\item{[5]}Dal Maso G., Skrypnik I.V.: Capacity theory for monotone operators.
{\it Potential Anal\/}. {\bf 7} (1997), 765-803.
\smallskip
\item{[6]}Dal Maso G., Skrypnik I.V.: Asymptotic behaviour of nonlinear
Dirichlet problems in perforated domains. {\it Ann. Mat. Pura Appl.\/} {\bf 
174} (1998), 13-72.
\smallskip
\item{[7]}Skrypnik I.V.: Nonlinear Elliptic Boundary Value Problems.
Teubner-Verlag, Leipzig, 1986.
\smallskip
\item{[8]}Skrypnik I.V.: Methods for Analysis of Nonlinear Elliptic
Boundary Value Problems. Nauka, Moscow, 1990. Transated in:
Translations of Mathematical Monographs 139, American Mathematical
Society, Providence, 1994.
\smallskip
\item{[9]}Skrypnik I.V.: New conditions for the homogenization of
nonlinear Dirichlet problems in perforated domains. {\it Ukrain. Mat. 
Zh.\/}, v.48, N 5, 1996.
\smallskip
\item{[10]}Skrypnik I.V.: Homogenization of nonlinear Dirichlet problems
in perforated domains of general structure. {\it Mat. Sb. (N.S.)\/} 
{\bf 187} (1996), 125-157.

\enddocument